
\vsize 24,5true cm
\vglue -1cm
\overfullrule 0mm
\parindent=0pt


\input amssym.def
\input amssym.tex


\def\makefootline{\baselineskip=52pt\line{\the\footline}}

\def\comp{\lbrack\kern-0.8mm\lbrack}

\newdimen\emm
\def\pmb#1{\emm=0.03em\leavevmode\setbox0=\hbox{#1}
\kern0.901\emm\raise0.434\emm\copy0\kern-\wd0
\kern-0.678\emm\raise0.975\emm\copy0\kern-\wd0
\kern-0.846\emm\raise0.782\emm\copy0\kern-\wd0
\kern-0.377\emm\raise-0.000\emm\copy0\kern-\wd0
\kern0.377\emm\raise-0.782\emm\copy0\kern-\wd0
\kern0.846\emm\raise-0.975\emm\copy0\kern-\wd0
\kern0.678\emm\raise-0.434\emm\copy0\kern-\wd0
\kern\wd0\kern-0.901\emm}

\newcount\coefftaille \newdimen\taille
\newdimen\htstrut \newdimen\wdstrut
\newdimen\ts \newdimen\tss

\font\caractsupgras=phvb
\def\Bbf#1{{\caractsupgras#1}}

\font\nomgot=eufm10 at 12pt

\def\fspeciale{\textfont0=\tenrmp%
\scriptfont0=\sevenrmp%
\scriptscriptfont0=\fivermp%
\textfont1=\tenip%
\scriptfont1=\sevenip%
\scriptscriptfont1=\fiveip%
\textfont2=\tensyp%
\scriptfont2=\sevensyp%
\scriptscriptfont2=\fivesyp%
\textfont3=\tenexp%
\scriptfont3=\tenexp%
\scriptscriptfont3=\tenexp%
\textfont\itfam=\tenitp%
\textfont\bffam=\tenbfp%
\textfont\slfam=\tenbfp%
\def\it{\fam\itfam\tenitp}%
\def\bf{\fam\bffam\tenbfp}%
\def\rm{\fam0\tenrmp}%
\def\sl{\fam\slfam\tenslp}%
\normalbaselineskip=12pt%
\multiply \normalbaselineskip by \coefftaille%
\divide \normalbaselineskip by 1000%
\normalbaselines%
\abovedisplayskip=10pt plus 2pt minus 7pt%
\multiply \abovedisplayskip by \coefftaille%
\divide \abovedisplayskip by 1000%
\belowdisplayskip=7pt plus 3pt minus 4pt%
\multiply \belowdisplayskip by \coefftaille%
\divide \belowdisplayskip by 1000%
\setbox\strutbox=\hbox{\vrule height\htstrut depth\wdstrut width 0pt}%
\rm}

\def\hfltoutpetit#1#2{{\mathop{\hbox to 0,5cm{\rightarrowfill}}
\limits^{\scriptstyle#1}_{\scriptstyle#2}}}

\def\fle{\rightarrow}

\def\vmid#1{\mid\!#1\!\mid}

\def\build#1_#2^#3{\mathrel{\mathop{\kern 0pt#1}\limits_{#2}^{#3}}}

\def\remonte{\vskip-\belowdisplayskip
\vskip-\baselineskip}

\def\messages#1{\immediate\write16{#1}}

\newdimen\margeg  \margeg=0pt
\def\findem{\vrule height0pt width4pt depth4pt}

\long\def\demApetit#1{{\parindent=0pt\messages{debut de preuve}\smallbreak
      \advance\margeg by 1truecm \leftskip=\margeg  plus 0pt
      {\everypar{\leftskip =\margeg  plus 0pt}
               \everydisplay{\displaywidth=\hsize
               \advance\displaywidth  by -1truecm
               \displayindent= 1truecm}
      {\bf Proof } -- \enspace #1
       \hfill\findem}\bigbreak}\messages{fin de preuve}}

\newdimen\itemgindent
\def\itemg#1{\itemgindent=\parindent
\advance\itemgindent by -0,5em
\item{\hbox to \the\itemgindent{#1\hfill}}}

\def\carreblanc{\vbox{\hrule
    \hbox{\vrule height 2mm\kern2mm\vrule}\hrule}}


\centerline{\Bbf{POSITIVE PARTITION RELATIONS FOR
\pmb{$P_\kappa(\lambda)$}}}
\bigskip
\centerline{\Bbf{PIERRE MATET and SAHARON SHELAH*}}\footnote{}{*This 
research was
supported by the Israel Science Foundation. Publication 804.

}

\vskip 1,5cm
{\bf Abstract.} \ Let $\kappa$ a regular uncountable cardinal and \ $\lambda$ \
a cardinal \ $>\kappa$, \ and suppose \ $\lambda^{<\kappa}$ \ is less than the
covering number for category \ $\hbox{{\bf cov}}(\pmb{M}_{\kappa,\kappa}).$ \
Then (a) $I_{\kappa,\lambda}^+ \hfltoutpetit{\kappa}{}
(I_{\kappa,\lambda}^+,\omega+1)^2,$ (b) $I_{\kappa,\lambda}^+
\hfltoutpetit{\kappa}{} [I_{\kappa,\lambda}^+]_{\kappa^+}^2$
if $\kappa$ is a limit cardinal, and (c) $I_{\kappa,\lambda}^+
\hfltoutpetit{\kappa}{} (I_{\kappa,\lambda}^+)^2$ if $\kappa$ is weakly
compact.
\vskip 1cm
{\Bbf{0. Introduction}}\footnote{}{2000 Mathematics Subject
Classification : 03E02, 03E35, 03E55, 03E04

Key words and phrases~:
$P_\kappa(\lambda),$ partition relation, weakly compact cardinal.}
\bigskip
Let $\kappa$ be a weakly compact cardinal. Then $\kappa\hfltoutpetit{}{}
(\kappa)^2$ and more generally for any cardinal
$\lambda\geq\kappa,\break
\{P_\kappa(\lambda)\}\hfltoutpetit{\kappa}{}
(I_{\kappa,\lambda}^+)^2$ ([M4]), which means that for any $F
: \kappa\times P_\kappa(\lambda) \hfltoutpetit{}{} 2,$ there is
$A\subseteq P_\kappa(\lambda)$ such that $A$ does not belong to
$I_{\kappa,\lambda}$ (the ideal of noncofinal subsets of
$P_\kappa(\lambda))$ and $F$ is constant on
$$\{(\cup(a\cap\kappa),b) : a,b\in A \hbox{ and }
\cup(a\cap\kappa)<\cup(b\cap\kappa)\}.$$
Now if $J$ is the ideal of noncofinal subsets of $\kappa,$
then $J^+ \hfltoutpetit{}{} (J^+)^2$ since $(A,<)$ is isomorphic to
$(\kappa,<)$ for any $A\in J^+.$ So it is natural to ask
whether $I_{\kappa,\lambda}^+\hfltoutpetit{\kappa}{}
(I_{\kappa,\lambda}^+)^2$ for every $\lambda>\kappa.$ It turns
out that the answer is negative. This is not surprising since
it is well-known that some members of $I_{\kappa,\lambda}^+$
may be quite different from $P_\kappa(\lambda).$ To give an
example , if the  \ GCH \ holds and $\lambda$ is the successor of a
cardinal of cofinality $<\kappa,$ then
$\overline{cof}(I_{\kappa,\lambda}\mid A)
<\overline{cof}(I_{\kappa,\lambda})$ for some $A\in
I_{\kappa,\lambda}^+$ ([MP\'eS2 ]). We prove that
$I_{\kappa,\lambda}^+\hfltoutpetit{\kappa}{}
(I_{\kappa,\lambda}^+)^2$ \ if and only if \ $\lambda^{<\kappa}$ \  is
less than \  $\hbox{{\bf cov}}(\pmb{M}_{\kappa,\kappa})$
(a generalization of the covering number for category $\hbox{{\bf
cov}}(\pmb{M})).$

\bigskip
Let $\kappa$ be an arbitrary regular uncountable cardinal.
Dushnik and Miller [DMi] established that
$\kappa\hfltoutpetit{}{}(\kappa,\omega)^2.$ This was improved to
$\kappa\hfltoutpetit{}{}(\kappa,\omega+1)^2$ by Erd\"os and Rado 
[ER]. The Erd\"os-Rado
result generalizes ([M3]) : for every cardinal \ $\lambda\geq \kappa,
\{P_\kappa(\lambda)\}\hfltoutpetit{\kappa}{} (I_{\kappa,\lambda}^+,\omega+1)^2$
\  (i.e. for any \ $F : \kappa\times 
P_\kappa(\lambda)\hfltoutpetit{}{} 2,$ \ there is
either
\ $A\in I_{\kappa,\lambda}^+$ \ such that \ $F$ \ is identically 0 on
$$\{(\cup(a\cap\kappa),b) : a,b\in A \hbox{ and }
\cup(a\cap\kappa)<\cup(b\cap\kappa)\},$$
or \ $a_0,a_1,\dots,a_\omega$ \ in \ $P_\kappa(\lambda)$ \ such that \
$a_0\subset a_1\subset\dots\subset a_\omega,
\cup(a_0\cap\kappa)<\cup(a_1\cap\kappa)<\dots<\cup(a_\omega\cap\kappa)$ \ and \
$F$ \ is identically 1 on \ $\{(\cup(a_n\cap\kappa),a_q) : n<q\leq\omega\}).$ \
Here we show that \ $I_{\kappa,\lambda}^+\hfltoutpetit{\kappa}{}
(I_{\kappa,\lambda}^+,\omega+1)^2$ \ if \ $\lambda^{<\kappa}$ \ is less than \
$\hbox{{\bf cov}}(\pmb{M}_{\kappa,\kappa}).$ \ In the other direction we prove
that \ $I_{\kappa,\lambda}^+ \hfltoutpetit{\kappa}{}\kern-3mm/\kern2mm
(I_{\kappa,\lambda}^+,3)^2$ \ if \ $\lambda$ \ is greater than or equal to \
${\hbox{\nomgot d}}_{\kappa}$ \ (or even \ $\overline{\hbox{\nomgot
d}}_{\kappa}).$
\bigskip
It is a result of [M5] that \
$\{P_\kappa(\lambda)\}\hfltoutpetit{\kappa}{}\kern-3mm/\kern2mm[I_{\kappa,\lambda}^+]_\lambda^2$
\ for any \ $\lambda>\kappa$ \ if \ $\kappa$ \ is a successor 
cardinal such that
\ $\kappa\hfltoutpetit{}{}\kern-3mm/\kern2mm[\kappa]_\kappa^2.$ \ In 
contrast to this,
we show that \ $I_{\kappa,\lambda}^+\hfltoutpetit{\kappa}{}
[I_{\kappa,\lambda}^+]_{\kappa^+}^2$ \ if \ $\kappa$ \ is a limit 
cardinal and \
$\lambda$ \ a cardinal \ $>\kappa$ \ with \ $\lambda^{<\kappa}< \hbox{{\bf
cov}}(\pmb{M}_{\kappa,\kappa}).$ \ It is also shown that \ 
$I^+_{\kappa,\lambda}
\hfltoutpetit{\kappa}{}\kern-3mm/\kern2mm[I^+_{\kappa,\lambda}]_\lambda^2$ 
\ if \
$\lambda\geq \overline{\hbox{\nomgot
d}}_{\kappa}.$
\bigskip
{\bf Throughout the remainder of this paper \pmb{$\kappa$} will denote a
regular uncountable cardinal and \pmb{$\lambda$} a cardinal
\pmb{$>\kappa.$}}
\bigskip
The paper is organized as follows. Section 1 reviews a number of
standard definitions concerning ideals on $\kappa$ and
$P_\kappa(\lambda).$ \ Sections 2-7 give results about combinatorics on \
$\kappa$ \ that are needed for our study of \ $P_\kappa(\lambda).$ \ Sections 2
and 3 review some facts concerning, respectively, the dominating number \
$\hbox{\nomgot d}_\kappa$ \ and the covering number for category \ $\hbox{{\bf
cov}}(\pmb{M}_{\kappa,\kappa}).$ \ Section 4 deals with the problem of
determining the value of the unequality number \ 
$\hbox{\nomgot{U}}_\kappa$ \ in the
case where \ $\kappa$ \ is a successor cardinal. In Section 5 we show that if \
$2^{<\kappa}=\kappa$ \ and \ $\hbox{\nomgot{U}}_\kappa 
<\kappa^{+\omega},$ \ then \
$\hbox{\nomgot{U}}_\kappa = \hbox{\pmb{non}}_\kappa$(weakly 
selective). Sections 6 and 7
review some material concerning, respectively, the unbalanced 
partition relation \
$J^+ \hfltoutpetit{}{} (J^+,\rho)^2$ \ and the square bracket 
partition relation \
$J^+\hfltoutpetit{}{} [J^+]_\rho^2.$
\bigskip
Sections 8-15 are concerned with combinatorial properties of ideals on \
$P_\kappa(\lambda).$ \ Section 8 gives two characterizations of $\hbox{\nomgot
d}_{\kappa,\lambda}^\kappa$ : one as the least cofinality of any
$\kappa$-complete fine ideal on $P_\kappa(\lambda)$ that is not a weak
$\pi$-point, and the other as the least cofinality of any
$\kappa$-complete fine ideal on $P_\kappa(\lambda)$ that admits a
maximal almost disjoint family of size $\kappa.$ \ In Section 9 we show
that any $\kappa$-complete fine ideal on
$P_\kappa(\lambda)$ with cofinality $<\hbox{{\bf
cov}}(\pmb{M}_{\kappa,\kappa})$ is a weak $\chi$-point. Conversely if 
\ $\kappa$
\ is inaccessible and \ $I_{\kappa,\lambda}$ \ is a weak \ $\chi$-point, then
$cof(I_{\kappa,\lambda})<\hbox{{\bf
cov}}(\pmb{M}_{\kappa,\kappa}).$ \ Sections 10-13 deal with 
unbalanced partition
relations. Given an infinite cardinal \ $\theta\leq\kappa$ \ such 
that \ $\kappa\fle
(\kappa,\theta)^2,$ \ we show that (a)
$u(\kappa,\lambda)\cdot\hbox{\pmb{non}}_\kappa(J^+\hfltoutpetit{}{}(J^+,\theta)^2)$ 
\ is the
least cofinality of any \ $\kappa$-complete fine ideal \ $H$ \ on \ 
$P_\kappa(\lambda)$ \
such that \ 
$H^+\hfltoutpetit{\kappa}{\kappa}\kern-3mm/\kern2mm(H^+;\theta)^2,$ \ 
(b) If
\ $H$ \ is a \ $\kappa$-complete fine ideal on \ $P_\kappa(\lambda)$ \ with \
$cof(H)<\hbox{{\bf cov}}(\pmb{M}_{\kappa,\kappa})$ \ (respectively,
$cof(H)<\hbox{\pmb{non}}_\kappa$(weakly selective)), then \
$H^+\hfltoutpetit{}{}(H^+,\theta)^2$ \ (respectively, \ 
$H^+\hfltoutpetit{\kappa}{\kappa}
(H^+,\theta)^2),$
\ and (c) Conversely, if \ $\theta=\kappa$ \ and \
$I_{\kappa,\lambda}^+\hfltoutpetit{\kappa}{\kappa}\ (H^+,\theta)^2,$ \ then \
$cof(I_{\kappa,\lambda})<\hbox{{\bf
cov}}(\pmb{M}_{\kappa,\kappa}).$ \ The last two sections are concerned with
square bracket partition relations. We show that if \ $\kappa$  \ is a limit
cardinal, then \ $H^+\hfltoutpetit{\kappa}{}[H^+]_{\kappa^+}^2$ \ 
(respectively,
\ $H^+\hfltoutpetit{\kappa}{\kappa}[H^+]_\kappa^2)$ \ for every ideal 
\ $H$ \ on \
$P_\kappa(\lambda)$ \ such that \ $cof(H)<\hbox{{\bf
cov}}(\pmb{M}_{\kappa,\kappa})$ \ (respectively, \
$cof(H)<\hbox{\pmb{non}}_\kappa(J^+\hfltoutpetit{}{}[J^+]_\kappa^2)).$ 
\ In the other
direction,
\  $\lambda\geq \overline{\hbox{\nomgot{d}}}_\kappa$ \ implies that \
$I_{\kappa,\lambda}^+
\hfltoutpetit{\kappa}{}\kern-3mm/\kern2mm[I_{\kappa,\lambda}^+]_\lambda^2$ \
(and \
$I_{\kappa,\lambda}^+\hfltoutpetit{\kappa}{\kappa}\kern-3mm/\kern2mm[I_{\kappa,\lambda}^+]_\kappa^2$
\ if \ $\kappa$ \ is a limit cardinal such that \ $2^{<\kappa} = \kappa).$
\vskip 1,5cm

{\Bbf{1. Ideals}}
\bigskip
In this section we review some standard definitions and a few basic
facts concerning ideals on $\kappa$ and $P_\kappa(\lambda).$
\bigskip
Given a cardinal $\mu$ and a set $A,$ let $P_\mu(A) = \{a\subseteq A :
\vmid a<\mu\}.$
\bigskip
Given an infinite set $S,$ an {\it ideal on} $S$ is a collection $K$
of subsets of $S$ such that (i) $\{s\}\in K$ for every $s\in S,$ (ii)
$P(A)\subseteq K$ for every $A\in K,$ (iii) $A\cup B\in K$ whenever
$A,B\in K,$ and (iv) $S\notin K.$
\bigskip
Given an ideal $K$ on $S,$ let $K^+ = P(S)-K$ and $K\mid A =
\{B\subseteq S : B\cap A\in K\}$ for $A\in K^+.$  $sat (K)$ is the
least cardinal $\tau$ with the property that for every $Y\subseteq K^+$
with $\vmid Y = \tau,$ there exist $A,B\in Y$ such that $A\not= B$ and
$A\cap B\in K^+.$
\bigskip
$cof(K)$ is the least cardinality of any $X\subseteq K$ such that $K =
\displaystyle\bigcup_{A\in X}P(A).\ K$ is {\it $\kappa$-complete} if
$\bigcup X\in K$ for every $X\in P_\kappa(K).$ Assuming that $K$ is
$\kappa$-complete and $\bigcup Y\in K^+$ for some $Y\subseteq K$ with
$\vmid Y = \kappa,\ \overline{cof}(K)$ is the least cardinality of any
$X\subseteq K$ such that $K = \bigcup \{P(\cup x) : x\in
P_\kappa(X)\}.$
\bigskip
{\bf We adopt the convention that the phrase ``ideal on
\pmb{$\kappa$}" means ``\pmb{$\kappa$}-complete ideal on
\pmb{$\kappa$}".}
\bigskip
Note that the smallest ideal on $\kappa$ is $P_\kappa(\kappa).$
\bigskip
Given two sets $A$ and $B$ and $f\in{}^AB,\ f$ is {\it regressive} if
$f(a)\in a$ for all $a\in A.$
\bigskip
An ideal $J$ on $\kappa$ is {\it normal} if given $A\in J^+$ and a
regressive $f\in{}^A\kappa,$ there is $B\in J^+\cap P(A)$ such that $f$
is constant on $B.$
\bigskip
$NS_\kappa$ \ denotes the nonstationary ideal on \ $\kappa.$
\bigskip
$\kappa$ \ is {\sl inaccessible} if \ $2^\mu<\kappa$ \ for every 
cardinal \ $\mu<\kappa.$
\bigskip
Let \ $[A]^2 = \{(\alpha,\beta)\in A\times A : \alpha<\beta\}$ \  for any
$A\subseteq\kappa.$ Given an ordinal $\alpha\geq 2,
\kappa\hfltoutpetit{}{}(\kappa,\alpha)^2$ \  means that for every \ $f :
[\kappa]^2\hfltoutpetit{}{} 2,$ \ there is \ $A\subseteq\kappa$ \ 
such that either \ $A$
\ has order type \ $\kappa$ \ and \ $f$ \ is identically 0 on
\ $[A]^2,$ \ or \ $A$ \ has order type \ $\alpha$ \ and \ $f$ \ is 
identically 1 on
\ $[A]^2.$ \ The negation of this and other partition relations is indicated by
crossing the arrow. \ $\kappa\hfltoutpetit{}{}(\kappa)^2$ \ means that
\ $\kappa\hfltoutpetit{}{}(\kappa,\kappa)^2.$
\bigskip
$\kappa$ is {\it weakly compact} if $\kappa\hfltoutpetit{}{}(\kappa)^2.$
\bigskip
If $\kappa$ is weakly compact, then it is inaccessible (see e.g. 
Proposition 4.4
in [Ka]).
\bigskip
An ideal \ $J$ \ on \ $\kappa$ \ is a {\sl weak} \ $P$-{\sl point} if 
given \ $A\in
J^+$ \ and \ $f\in{}^A\kappa$ \ with \ $\{f^{-1}(\{\alpha\}) :
\alpha<\kappa\}\subseteq J,$ \ there is \ $B\in J^+\cap P(A)$ \ such 
that \ $f$ \ is
\ $<\kappa$-to-one on \ $B.$ \ $J$ \ is a {\sl local} \ $Q$-{\sl 
point} if given \
$g\in{}^\kappa\kappa,$ \ there is \ $B\in J^+$ \ such that \ 
$g(\alpha)<\beta$ \ for
any \ $(\alpha,\beta)\in[B]^2.$ \ $J$ \ is a {\sl weak} \ $Q$-{\sl 
point} if \ $J\mid
A$ is a local \ $Q$-point for every \ $A\in J^+.$
\bigskip
It is well-known (see [M1] for a proof) that an ideal \ $J$ \ on \ 
$\kappa$ \ is a
weak \ $Q$-point if and only if given \ $A\in J^+$ \ and a \ 
$<\kappa$-to-one \ $f :
A\hfltoutpetit{}{}\kappa,$ \ there is \ $B\in J^+\cap P(A)$ \ such that 
\ $f$ \ is
one-to-one on \
$B.$
\bigskip
An ideal \ $J$ \ on \ $\kappa$ \ is {\sl weakly selective} if it is 
both a weak \
$P$-point and a weak \ $Q$-point.
\bigskip
Given a cardinal \ $\rho$ \ with \ $2\leq\rho\leq\kappa$ \ and an 
ideal \ $J$ \ on \
$\kappa, J^+\hfltoutpetit{}{}[J^+]_\rho^2$ \ means that for every \ 
$A\in J^+$ \ and
every \ $f : [A]^2\hfltoutpetit{}{}\rho,$ \ there is \ $B\in J^+\cap 
P(A)$ \ such that \
$f''[B]^2\not=\rho. \ \kappa\hfltoutpetit{}{}[\kappa]_\rho^2$ \ means that \
$(P_\kappa(\kappa))^+\hfltoutpetit{}{}[(P_\kappa(\kappa))^+]_\rho^2.$
\bigskip
  Note that \ $\kappa\hfltoutpetit{}{}[\kappa]_2^2$ \ if and only if \
$\kappa\hfltoutpetit{}{}(\kappa)^2.$
\bigskip
Let \ $P$ \ be a property such that at least one ideal on \ $\kappa$ \ does not
satisfy \ $P.$ \ Then \ $\hbox{\bf non}_\kappa(P)$ \ (respectively,
$\overline{\hbox{\bf non}}_\kappa(P))$
\ denotes the least cardinal \ $\tau$ \ for which one can find an 
ideal \ $J$ \ on \
$\kappa$ \ such that \ $cof(J) = \tau$ \ (respectively, \ 
$\overline{cof}(J) = \tau)$
\ and  \ $J$ \ does not satisy \ $P.$
\bigskip
Notice that \ 
$\lambda^{<\kappa}<\overline{\hbox{\pmb{non}}}_\kappa(P)$ \ if and 
only if \
$\lambda^{<\kappa}<\hbox{\pmb{non}}_\kappa(P).$
\bigskip
$I_{\kappa,\lambda}$ denotes the set of all $A\subseteq
P_\kappa(\lambda)$ such that $A\cap\{b\in P_\kappa(\lambda) :
a\subseteq b\} = \phi$ for some $a\in P_\kappa(\lambda).$ An ideal $H$
on $P_\kappa(\lambda)$ is {\it fine} if $I_{\kappa,\lambda}\subseteq H.$
\bigskip
{\bf We adopt the convention that the phrase ``ideal on
\pmb{$P_\kappa(\lambda)$}" means ``\pmb{$\kappa$}-complete fine ideal
on \pmb{$P_\kappa(\lambda)$}".}
\bigskip
Note that $I_{\kappa,\lambda}$ is the smallest ideal on
$P_\kappa(\lambda).$
\bigskip
$u(\kappa,\lambda)$ denotes the least cardinality of any $A\in
I_{\kappa,\lambda}^+.$
\bigskip
The following facts are well-known (see e.g.
[MP\'eS1]) : (1) $u(\kappa,\lambda)\geq\lambda$ ; (2) $\lambda^{<\kappa} =
2^{<\kappa}\cdot u(\kappa,\lambda)$ ; (3) $u(\kappa,\lambda) =
cof(I_{\kappa,\lambda}\mid A)$ for every $A\in I_{\kappa,\lambda}^+$ ;
(4) $u(\kappa,\kappa^{+n}) = \kappa^{+n}$ whenever $0<n<\omega.$
\bigskip
${\cal K}(\kappa,\lambda)$ \ denotes the set of all cardinals \
$\sigma\geq\lambda$ \ with the property that there is \ $T\subseteq
P_\kappa(\lambda)$ \ such that \ $\vmid T = \sigma$ \ and \ 
$\vmid{T\cap P(a)}<\kappa$
\ for every \ $a\in P_\kappa(\lambda).$
\bigskip
It is simple to see that \ $\sigma\leq u(\kappa,\lambda)$ \ for every \
$\sigma\in{\cal K}(\kappa,\lambda).$ \ Notice that \ $\lambda\in{\cal
K}(\kappa,\lambda).$ \ More generally, if \ $\tau$ \ is an infinite cardinal \
$\leq\kappa$ \ such that \ $\vmid{P_\tau(\nu)}<\kappa$ \ for every 
infinite cardinal \
$\nu<\kappa,$ \ then \ $\lambda^{<\tau}\in{\cal K}(\kappa,\lambda).$ \ It
follows that \ $\lambda^{<\kappa}\in{\cal K}(\kappa,\lambda)$ \ if \ $\kappa$
\ is inaccessible. It can be shown (see Remark 11.4  in [To 2] and 
Theorem 4.1 in
[CFMag]) that \ $\lambda^+\in{\cal K}(\kappa,\lambda)$ \ if \
$\carreblanc_\kappa^*$
\ holds and \ $cf(\lambda)<\kappa.$
\bigskip
An ideal $H$ on $P_\kappa(\lambda)$ is {\it $\kappa$-normal} if given
$A\in H^+$ and a regressive $f\in{}^A\kappa,$ there is $B\in H^+\cap
P(A)$ such that $f$ is constant on $B.$ The smallest $\kappa$-normal
ideal on $P_\kappa(\lambda)$ is denoted by $NS_{\kappa,\lambda}^\kappa.$
\vskip 1,5cm
{\Bbf{2. Domination}}
\bigskip
In this section we recall some characterizations of the dominating number \
$\hbox{\nomgot d}_\kappa.$
\bigskip
{\bf Definition.} \ $\hbox{\nomgot d}_\kappa$ \ is the least 
cardinality of any \
$X\subseteq {}^\kappa\kappa$ \ with the property that for every \
$g\in{}^\kappa\kappa,$ \ there is \ $f\in X$ \ such that \ 
$g(\alpha)<f(\alpha)$ \
for all \ $\alpha<\kappa.$

$\overline{\hbox{\nomgot d}}_\kappa$ \ is the least cardinality of 
any \ $X\subseteq
{}^\kappa\kappa$ \ with the property that for every \ $g\in 
{}^\kappa\kappa,$ \ there
is \ $x\in P_\kappa(X)$ \ such that \ $\displaystyle 
g(\alpha)<\bigcup_{f\in x}f(\alpha)$ \
for all
\ $\alpha<\kappa.$
\bigskip\bigskip
{\bf PROPOSITION 2.1.}
\medskip
{\parindent=1,5cm
\item{\bf (i)} ([L1]) \ $\hbox{\nomgot d}_\kappa = cof(NS_\kappa).$
\smallskip
\item{\bf (ii)} ([MRoS]) \ $\overline{\hbox{\nomgot d}}_\kappa =
\overline{cof}(NS_\kappa).$
\par}
\bigskip\bigskip
{\bf Definition.} \ Given an ideal \ $J$ \ on \ $\kappa,\ {\cal
M}_J^{\geq\kappa}$ \ is the set of all \ $Q\subseteq J^+$ \ such that 
(i) $\vmid
Q\geq\kappa,$ (ii) $A\cap B\in J$ \ for all \ $A,B\in Q$ \ with \ 
$A\not= B,$ \ and
(iii) for every \ $C\in J^+,$ \ there is \ $A\in Q$ \ with \ $A\cap C\in J^+.$

$\hbox{\nomgot a}_J$ \ is the least cardinality of any member of \ ${\cal
M}_J^{\geq\kappa}$ \ if \ ${\cal M}_J^{\geq\kappa}\not=\phi,$ \ and \
$(2^\kappa)^+$ \ otherwise.
\bigskip\bigskip
{\bf THEOREM 2.2.} ([Laf], [MP2]) \ $\hbox{\nomgot d}_\kappa =
\hbox{\pmb{non}}_\kappa(\hbox{\nomgot a}_J>\kappa) = 
\hbox{\pmb{non}}_\kappa$(weak \
$P$-point).
\bigskip\bigskip
{\bf PROPOSITION 2.3.} \   $\overline{\hbox{\nomgot
d}}\kappa\geq\overline{\hbox{\pmb{non}}}_\kappa(\hbox{\nomgot
a}_J>\kappa)\geq\overline{\hbox{\pmb{ non}}}_\kappa$(weak \ $P$-point).
\bigskip
{\bf Proof.} \ The first inequality follows from Proposition 2.1 (ii)  since \
$\hbox{\nomgot a}_{{NS}_\kappa} =
\kappa$ ([MP2]). To prove the second inequality, argue as for Lemma 8.5 below.
\hfill$\carreblanc$
\bigskip\bigskip
{\bf QUESTION.} \ Is it consistent that \ $\overline{\hbox{\nomgot
d}}_\kappa>\overline{\hbox{\pmb{non}}}_\kappa$(weak \ $P$-point) ?
\vskip 1,5cm
{\Bbf{3. Covering for category}}
\bigskip
Throughout this section \ $\nu$ \ will denote  a fixed regular 
infinite cardinal.
\bigskip
We will review some basic facts concerning the covering number \ $\hbox{{\bf
cov}}(\pmb{M}_{\nu,\nu}).$
\bigskip\bigskip
{\bf Definition.} \ Suppose \ $\rho$ \ is a cardinal \ $\geq\nu.$

Let \ $Fn(\rho,2,\nu) = \cup\{{}^a2 : a\in P_\nu(\rho)\}. \ 
Fn(\rho,2,\nu)$ \ is ordered
by~: $p\leq q$ \ if and only  if \ $q\subseteq p.$

${}^\rho2$ \ is endowed with the topology obtained by taking as basic 
open sets \
$\phi$ \ and \ $O_s^\rho$ \ for \ $s\in Fn(\rho,2,\nu),$ \ where \ $O_s^\rho =
\{f\in{}^\rho 2 : s\subseteq f\}.$

$\pmb{M}_{\nu,\rho}$ \ is the set of all \ $W\subseteq{}^\rho2$ \ such that \
$W\cap(\cap X) = \phi$ \ for some collection \ $X$ \ of dense open 
subsets of \ ${}^\rho2$
\ with \ $0<\vmid X\leq\nu.$

$\hbox{{\bf cov}}(\pmb{M}_{\nu,\rho})$ \ is the least cardinality of 
any \ $Y\subseteq
\pmb{M}_{\nu,\rho}$ \ such that \ ${}^\rho2 = \cup Y.$
\bigskip\bigskip
{\bf PROPOSITION 3.1.}
\medskip
{\parindent=1,5cm
\item{\bf (i)} ([L2],[Mil2]) \   $\hbox{{\bf
cov}}(\pmb{M}_{\nu,\rho})\geq\nu^+$ \ {\sl for every cardinal \ $\rho\geq\nu.$}
\medskip
\item{\bf (ii)} ([L2],[Mil2]) \  {\sl Suppose that \ $\rho$ \ and \ $\mu$ \ are
two cardinals such that \ $\nu\leq\mu\leq\rho.$ \ Then }\ $\hbox{{\bf
cov}}(\pmb{M}_{\nu,\mu})\geq \hbox{{\bf cov}}(\pmb{M}_{\nu,\rho}).$
\medskip
\item{\bf (iii)} ([L2]) \  {\sl Suppose \ $2^{<\nu}>\nu.$ \ Then }\ $\hbox{{\bf
cov}}(\pmb{M}_{\nu,\nu}) = \nu^+.$
\par}
\bigskip\bigskip
{\bf PROPOSITION 3.2.} \ {\sl Suppose that \ $\rho$ \ is a cardinal \ $>\nu$ \
and
\
$V\models 2^{<\nu} = \nu.$ \ Then setting \ $P = Fn(\rho,2,\nu)~:$}
\medskip
{\parindent=1,5cm
\item{\bf (i)}([L2],[Mil2]) \   $V^P\models\hbox{{\bf
cov}}(\pmb{M}_{\nu,\rho})\geq\rho.$
\medskip
\item{\bf (ii)} ([L2],[Mil2]) \  {\sl If \ $cf(\rho)\leq\nu,$ \ then}
$V^P\models\hbox{{\bf cov}}(\pmb{M}_{\nu,\nu})>\rho.$
\item{\bf (iii)}  {\sl Let \ $\mu$ \ be any regular cardinal \ $>\nu.$ \ Then \
$(\hbox{\nomgot d}_\mu)^{V^P} = (\hbox{\nomgot d}_\mu)^V$ \ and \ $
(\overline{\hbox{\nomgot d}}_\mu)^{V^P}\leq(\overline{\hbox{\nomgot 
d}}_\mu)^V.$}
\par}
\bigskip
{\bf Proof.} \ (iii) : The conclusion easily follows from the following
observation~: Suppose that \ $\sigma$ \ is a cardinal \ $>0$ \ and \ 
$F\in V^P$ \
is a function from \ $\sigma\times\mu$ \ to \ $\mu.$  \ Then by Lemma 
VII.6.8 of
\ [K], \ there is \ $H : \sigma\times\mu\hfltoutpetit{}{} 
P_{\nu^+}(\mu)$ \ such that \
$H\in V$ \ and \ $F(\alpha,\beta)\in H(\alpha,\beta)$ \ (so \
$F(\alpha,\beta)\leq\cup H(\alpha,\beta))$ \ for every \
$(\alpha,\beta)\in\sigma\times\mu.$ \hfill$\carreblanc$
\bigskip
{\bf Remark.} \ It is not known whether it is consistent that \ $cf({\hbox{\bf
cov}}(\pmb{M}_{\nu,\nu})\leq\nu.$
\vskip 1,5cm
{\Bbf{4. Unequality}}
\bigskip
Our main concern in this section is with the problem of evaluating 
the unequality
number \ $\hbox{\nomgot U}_\kappa$ \ when \ $\kappa$ \ is a successor cardinal.
\bigskip
{\bf Definition.} \ $\hbox{\nomgot U}_\kappa$ \ (respectively, \ $\hbox{\nomgot
U}'_\kappa)$ \ is the least cardinality of any \ $F\subseteq {}^\kappa\kappa$ \
with the property that for every \ $g\in{}^\kappa\kappa,$ \ there is 
\ $f\in F$ \
such that \ $\{\alpha\in\kappa : f(\alpha) = g(\alpha)\}=\phi$ \ 
(respectively \
$\vmid{\{\alpha\in\kappa : f(\alpha) = g(\alpha)\}}<\kappa).$
\bigskip
The following is readily checked.
\bigskip\bigskip
{\bf PROPOSITION 4.1.} \ $\hbox{\bf
cov}(\pmb{M}_{\kappa,\kappa})\leq\hbox{\nomgot 
U}_\kappa\leq\hbox{\nomgot d}_\kappa.$
\bigskip
{\bf Remark.} \ It is shown in [MRoS] that if \ $V\models$ GCH, \ then
there is a \ $\kappa$-complete \ $\kappa^+-$cc \ forcing notion \ $P$ 
\ such that
$$V^P\models ``\overline{\hbox{\nomgot d}}_\kappa = \kappa^{+\omega} 
\hbox{ and }
\hbox{\bf cov}(\pmb{M}_{\kappa,\kappa}) = 2^\kappa = \kappa^{+(\omega+1)}".$$
For models where \ $\hbox{\nomgot d}_\kappa>\kappa^+$ \ see also [CS].
\bigskip\bigskip
{\bf PROPOSITION 4.2.} \ $\hbox{\nomgot U}_\kappa = \hbox{\nomgot U}'_\kappa.$
\bigskip
{\bf Proof.} \ Fix $F\subseteq{}^\kappa\kappa$ with the property that for every
$g\in{}^\kappa\kappa,$ there is $f\in F$ such that
$$\vmid{\{\alpha\in\kappa :
f(\alpha) = g(\alpha)\}}<\kappa.$$
  For $f\in F$ and $\gamma,\delta<\kappa,$
define $f_{\gamma,\delta}\in{}^\kappa\kappa$ by~: $f_{\gamma,\delta} (\alpha) =
f(\alpha)$ if $\alpha\geq\gamma,$ and $f_{\gamma,\delta}(\alpha) = \delta$
otherwise. Then for every $g\in{}^\kappa\kappa,$ there are $f\in F$ and
$\gamma,\delta<\kappa$ such that $\{\alpha\in\kappa~:
f_{\gamma,\delta}(\alpha) = g(\alpha)\} = \phi.$ \hfill\carreblanc
\bigskip
The following is due to Landver [L2].
\bigskip\bigskip
{\bf PROPOSITION 4.3.}\ {\sl $cf(\hbox{\nomgot U}_\kappa)>\kappa.$}
\bigskip
{\bf Proof.} \ Suppose otherwise. Set $\nu = cf(\hbox{\nomgot 
U}_\kappa)$ and fix
$F\subseteq{}^\kappa\kappa$ so that $\vmid F = \hbox{\nomgot 
U}_\kappa$ and for every
$g\in{}^\kappa\kappa,$ there exists $f\in F$ with $\{\alpha\in\kappa :
f(\alpha)=g(\alpha)\}=\phi.$ Let
$<F_\beta :
\beta<\nu>$ be such that (a) $\vmid{F_\beta}<\hbox{\nomgot U}_\kappa$ 
for any $\beta,$ and
(b)
$\displaystyle\bigcup_{\beta<\nu}F_\beta = F.$ Select $A_\beta\subseteq\kappa$
for $\beta<\nu$ so that (i) $\vmid{A_\beta} = \kappa$ for every $\beta<\nu,$
(ii) $A_\beta\cap A_\gamma=\phi$ whenever $\gamma<\beta<\nu,$ and (iii)
$\displaystyle\bigcup_{\beta<\nu}A_\beta = \kappa.$ For each $\beta<\nu,$ there
is $g_\beta : A_\beta\hfltoutpetit{}{}\kappa$ such that
$$\{\alpha\in A_\beta : (f\upharpoonright A_\beta)(\alpha) =
g_\beta(\alpha)\}\not=\phi$$
for every $f\in F_\beta.$ Set $g = \displaystyle\bigcup_{\beta<\nu}g_\beta.$
Then clearly, $\{\alpha\in\kappa : f(\alpha) = g(\alpha)\}\not=\phi$ for all
$f\in F.$ This is a contradiction.\hfill$\carreblanc$\eject
\bigskip
We now turn our attention to the task of computing \ $\hbox{\nomgot 
U}_\kappa.$ \ We
begin with the case when \ $\kappa$ \ is a successor cardinal.
\bigskip\bigskip
{\bf THEOREM 4.4.} \ {\sl Suppose \ $\kappa$ \ is the successor of a 
regular infinite
cardinal \ $\nu.$ \ Then }
$$\hbox{\nomgot U}_\kappa \geq \hbox{min}(\hbox{\nomgot d}_\kappa,\hbox{\bf
cov}(\pmb{M}_{\nu,\kappa})).$$
\bigskip
{\bf Proof.} \ Fix \ $F\subseteq{}^\kappa\kappa$ \ with \
$0<\vmid F<\hbox{min}(\hbox{\nomgot d}_\kappa,\hbox{\bf 
cov}(\pmb{M}_{\nu,\kappa})).$ \
Pick \ $k : \kappa\hfltoutpetit{}{}\kappa-\nu$ \ so that
$$\vmid{\{\alpha<\kappa : k(\alpha)>f(\alpha)\}} = \kappa$$
for every \ $f\in F.$ \ Select a bijection \ $j : 
\kappa\times\nu\hfltoutpetit{}{}\kappa$
\ and a bijection \ $i_\alpha : k(\alpha)\hfltoutpetit{}{}\nu$ \ for each \
$\alpha<\kappa.$ \ Given \
$A\subseteq\kappa$ \ and \ $t\in{}^A2,$ \ define a partial function \ 
$\overline t$ \
from \ $\kappa$ \ to \ $\kappa$ \ by stipulating that \ $\overline 
t(\alpha) = \gamma$ \
if and only if (a) \ $\gamma<k(\alpha),$
(b) $\{j(\alpha,\eta) : \eta<i_\alpha(\gamma)\}\subseteq t^{-1}(\{0\}),$ \
and (c) \ $j(\alpha,i_\alpha(\gamma))\in t^{-1}(\{1\}).$ \  For \ 
$f\in F,$ \ let \
$D_f$ \ be the set of all \ $s\in Fn(\kappa,2,\nu)$ \ such that there 
is \ $\alpha\in
dom(\overline s)$ \ with \ $k(\alpha)>f(\alpha)$ \ and \ $\overline 
s(\alpha)=f(\alpha).$ \
Clearly, each \ $D_f$ \ is a dense subset of \ $Fn(\kappa,2,\nu),$ \ 
so we can find \
$g\in {}^\kappa2$ \ with the property that for every \ $f\in  F,$ \
there is\ $a\in P_\nu(\kappa)$ \ with \ $g\upharpoonright a\in D_f.$
\ Then
$$\{\alpha\in dom(\overline g) : \overline g(\alpha) =
f(\alpha)\}\not=\phi$$
for every \ $f\in F.$ \hfill\carreblanc
\bigskip\bigskip
{\bf THEOREM 4.5.} \ {\sl Suppose \ $\kappa$ \ is a successor cardinal.
Then }\ $\hbox{\nomgot U}_\kappa\geq\overline{\hbox{\nomgot
d}}_\kappa.$
\bigskip
{\bf Proof.} \ Fix \ $F\subseteq{}^\kappa\kappa$ \ with \ $0<\vmid
F<\overline{\hbox{\nomgot d}}_\kappa.$ \ Set \ $\kappa=\nu^+.$ \ Pick
\ $k : \kappa\hfltoutpetit{}{}\kappa-\nu$ \ so that
$$\vmid{\{\alpha<\kappa : f(\alpha)<k(\alpha)\}} = \kappa$$
for every \ $f\in F.$ \ For \ $\alpha<\kappa,$ \ select a bijection
\ $\pi_\alpha : k(\alpha) \hfltoutpetit{}{}\nu.$ \ Given \ $f\in F,$ \ there
exists \ $i_f\in\nu$ \ such that the set
$$A_f = \{\alpha<\kappa : f(\alpha)<k(\alpha) \hbox{ and }
\pi_\alpha(f(\alpha)) = i_f\}$$
has size \ $\kappa.$ \ Define \ $g_f\in{}^\kappa\kappa$ \ by
$$g_f(\beta) = \hbox{ least } \alpha\in A_f \hbox{ such that }
\alpha\geq\beta.$$
It  is shown in [MRoS] that \ $\overline{\hbox{\nomgot d}}_\kappa$ \
is the least cardinality of any \ $X\subseteq{}^\kappa\kappa$ \ with
the property that for every \ $h\in{}^\kappa\kappa,$ \ there is \
$x\in P_\kappa(X)$ \ such that the set
\ $\{\beta<\kappa : h(\beta)\geq\displaystyle\bigcup_{f\in x}f(\beta)\}$
is nonstationary in \ $\kappa.$ \ Hence there is \
$h\in{}^\kappa\kappa$ \ such that the set
$$B_x = \{\beta<\kappa : h(\beta)\geq\bigcup_{f\in x}g_f(\beta)\}$$
is stationary in \ $\kappa$ \ for every \ $x\in P_\kappa(F).$
\bigskip
Define \ $J\subseteq P(\kappa)$ \ by : $D\in J$ \ if and only if
there is \ $x\in P_\kappa(F)$ \ such that \ $D\cap B_x\in
NS_\kappa.$ \ Then \ $J$ \ is an ideal on \ $\kappa.$ \ Since \
$sat(J)>\nu$ \ by a result of Ulam (see [Ka], 16.3), there exist
pairwise disjoint \ $D_i\in J^+$ \ for \ $i<\nu$ \ with \
$\displaystyle\bigcup_{i<\nu}D_i = \kappa.$
\bigskip
Let \ $C$ \ be the set of all infinite limit ordinals \
$\delta<\kappa$ \ such that \ $h(\xi)<\delta$ \ for every \
$\xi<\delta.$ \ Then \ $C$ \ is a closed unbounded subset of \
$\kappa.$ \ Define \ $t\in{}^\kappa\kappa$ \ so that for every \
$\eta<\kappa, t(\eta)<k(\eta)$ \ and \ $c_\eta\in
D_{\pi_\eta(t(\eta))},$ \ where \ $c_\eta = \cup(C\cap\eta).$
\bigskip
Now fix \ $f\in F.$ \ Pick \ $\zeta\in D_{i_f}\cap C\cap B_{\{f\}}$ \
and set \ $\eta = g_f(\zeta).$ \ Notice that \ $\zeta\leq\eta$ \ by
the definition of \ $g_f.$ \ Also, \ $\eta\leq h(\zeta)$ \ since \
$\zeta\in B_{\{f\}}.$ \ Hence \ $c_\eta = \zeta$ \ by the definition
of \ $C$ \ and the fact that \ $\zeta\in C.$ \ It now follows from
the definition of \ $t$ \ and the fact that \ $\zeta\in D_{i_f}$ \
that \ $\pi_\eta(t(\eta)) = i_f.$ \ On the other hand, \ $\eta\in
A_f$ \ since \ $\eta = g_f(\zeta),$ \ so \ $f(\eta)<k(\eta)$ \ and \
$\pi_\eta(f(\eta)) = i_f.$ \ Thus \ $t(\eta) = f(\eta).$
\hfill\carreblanc
\bigskip
{\bf Remark.} \ It follows from Proposition 4.1 and Theorem 4.5 that
\ $\hbox{\nomgot U}_\kappa = \hbox{\nomgot d}_\kappa$ \ if \
$\kappa$ \ is a successor cardinal and \ $\hbox{\nomgot
d}_\kappa<\kappa^{+\omega}.$
\bigskip\bigskip
{\bf THEOREM 4.6.} \ {\sl Suppose that \ $\kappa$ \ is a successor
cardinal and \ $2^{<\kappa}=\kappa.$ \ Then \ $\hbox{\nomgot
U}_\kappa = \hbox{\nomgot d}_\kappa.$}
\bigskip
{\bf Proof.} \ By Proposition 4.1 it suffices to prove that \
$\hbox{\nomgot U}_\kappa\geq\hbox{\nomgot d}_\kappa.$ \ Set \
$\kappa=\nu^+$ \ and select a one-to-one
$$j : \bigcup_{\alpha<\kappa}
{}^{[\alpha,\alpha+\nu)}\kappa\hfltoutpetit{}{}\kappa.$$
Now fix \ $F\subseteq{}^\kappa\kappa$ \ with \ $0<\vmid
F<\hbox{\nomgot d}_\kappa.$ \
Select \ $g\in{}^\kappa\kappa$ \ so that for every \ $f\in F,$ \
there is \ $\beta_f<\kappa$ \ with
$$j\bigl(f\upharpoonright[\beta_f,\beta_f+\nu)) < g(\beta_f).$$
Let \ $C$ \ be the set of all \ $\gamma<\kappa$ \ such that \
$\beta+\nu<\gamma$ \ and \ $g(\beta)<\gamma$ \ for every \
$\beta<\gamma.$ \ Then \ $C$ \ is a closed unbounded subset of \
$\kappa.$ \ Let \ $<\gamma_\delta : \delta<\kappa>$ \ be the
increasing enumeration of \ $C.$ \ For \ $\delta<\kappa,$ \ set
$$W_\delta =
\Bigl\{t\in\bigcup_{\gamma_\delta\leq\alpha<\gamma_{\delta+1}}
{}^{[\alpha,\alpha+\nu)}\kappa:j(t)<\gamma_{\delta+1}\Bigr\}.$$
Then define \
$k_\delta\in{}^{[\gamma_\delta,\gamma_{\delta+1})}\kappa$ \ so
that for every \ $t\in W_\delta,$ \ there is \ $\zeta\in dom(t)$ \
with
\ $k_\delta(\zeta) = t(\zeta).$ \ Set \ $k =
\bigcup_{\delta<\kappa}k_\delta.$
\bigskip
Given \ $f\in F,$ \ let \ $\delta_f<\kappa$ \ be such that \
$\gamma_{\delta_f}\leq\beta_f<\gamma_{\delta_f+1}.$ \ Then \
$f\upharpoonright[\beta_f,\beta_f+\nu)\in W_{\delta_f}.$ \ Hence \
$k(\zeta) = f(\zeta)$ \ for some \
$\zeta\in[\beta_f,\beta_f+\nu).$ \hfill\carreblanc
\bigskip
{\bf QUESTION.} \ Is it consistent that \ $\kappa$ \ is a successor
cardinal and \ $\hbox{\nomgot U}_\kappa<\hbox{\nomgot d}_\kappa$ ?
\bigskip
{\bf QUESTION.} \ Is it consistent that \ $\kappa$ \ is a successor
cardinal such that \ $2^{<\kappa}=\kappa$ \ and \ $\hbox{\bf
cov}(\pmb{M}_{\kappa,\kappa}) < \hbox{\nomgot U}_\kappa$ ?
\bigskip
Let us now consider the case when \ $\kappa$ \ is a limit cardinal.
By a result of Bartoszy\'nski [B] and Miller [Mil1], \ $\hbox{\nomgot
U}_\omega = \hbox{\bf cov}(\pmb{M}_{\omega,\omega}).$ \ Landver [L2]
was able to show that this fact generalizes to uncountable
inaccessible cardinals~:
\bigskip\bigskip
{\bf THEOREM 4.7.} \ {\sl If \ $\kappa$ \ is an inaccessible cardinal,
then} \ ${\hbox{\nomgot U}}_\kappa = \hbox{\bf
cov}(\pmb{M}_{\kappa,\kappa}).$
\bigskip
{\bf QUESTION.} \ Is it consistent that \ $\kappa$ \ is a limit
cardinal and \ $\hbox{\bf cov}(\pmb{M}_{\kappa,\kappa})<\hbox{\nomgot
U}_\kappa$ ?
\vskip 1,5cm
{\Bbf{5. Weak selectivity}}
\bigskip
The following is due to Baumgartner, Taylor and Wagon [BauTW].
\bigskip
{\bf PROPOSITION 5.1.} \ {\sl If \ $\kappa$ \ is a successor
cardinal, then every ideal on \ $\kappa$ \ is a weak \ $Q$-point.}
\bigskip
By Proposition 5.1 and Theorem 2.2 \ $\hbox{\pmb{non}}_\kappa$(weakly
selective) $= \hbox{\nomgot d}_\kappa$ \ if \ $\kappa$ \ is a
successor cardinal. The remainder of the section is primarily
concerned with the value of \ $\hbox{\pmb{non}}_\kappa$(weakly 
selective) in the
case when \ $\kappa$ \ is a limit cardinal.
\bigskip
{\bf Remark.} \ It is easy to see that \ $\kappa^+\leq 
\hbox{\pmb{non}}_\kappa$(weak \
$Q$-point) if \ $\kappa$ \ is a limit cardinal.
\bigskip
{\bf Definition.} \ An ideal \ $J$ \ on \ $\kappa$ \ is a weak
semi-$Q$-point if given \ $A\in J^+$ \ and a \ $<\kappa$-to-one
function \ $f$ \ from \ $A$ \ to \ $\kappa,$ \ there is \ $C\in
J^+\cap P(A)$ \ such that \ $\vmid{C\cap
f^{-1}(\{\alpha\})}\leq\vmid{\alpha}$ \ for every \
$\alpha\in\kappa.$

$J$ \ is weakly semiselective if \ $J$ \ is both
a weak semi-$Q$-point and a weak \ $P$-point.

$J$ \ is weakly rapid if given \ $A\in J^+$ \ and \
$f\in{}^\kappa\kappa,$ \ there is \ $C\in J^+\cap P(A)$ \ such that
\ o.t.$(C\cap f(\alpha))\leq\alpha+1$ \ for every \
$\alpha\in\kappa.$
\bigskip
{\bf Remark.} \ It is simple to see that every weak \ $Q$-point
ideal on \ $\kappa$ \ is weakly rapid, and every weakly rapid ideal
on \ $\kappa$ \ is a weak semi-$Q$-point.
\bigskip
Every weak semi-$Q$-point ideal on \ $\omega$ \ is weakly rapid
([MP1]). We will show that this does not generalize.
\bigskip
{\bf Definition.} \ An ideal \ $J$ \ on \ $\kappa$ \ is a
semi-$Q$-point if given a \ $<\kappa$-to-one function \ $f$ \ from \
$\kappa$ \ to \ $\kappa,$ \ there is \ $B\in J$ \ such that \
$\vmid{f^{-1}(\{\alpha\})-B}\leq\vmid\alpha$ \ for every \
$\alpha\in\kappa.$
\bigskip\bigskip
{\bf PROPOSITION 5.2.} \ {\sl Suppose \ $\kappa$ \ is a limit cardinal.
Then there exists a semi-$Q$-point ideal on \ $\kappa$ \ that is not
weakly rapid.}
\bigskip
{\bf Proof.} \ Let \ $Y$ \ be the set of all infinite cardinals \
$<\kappa.$ \ Select \ $h\in{}^Y\kappa$ \ so that (a) $h(\mu)$ \ is a
regular infinite cardinal \ $\leq\mu$ \ for every \ $\mu\in Y,$ \
and (b) $\{\mu\in Y$~: $h(\mu)\geq\theta\}$ \ is stationary in \
$\kappa$ \ for every \ $\theta\in Y.$ \ For \ $A\subseteq\kappa$ \
and \ $\theta\in Y,$ \ let \ $T_\theta^A$ \ be the set of all \
$\mu\in Y$ \ such that \ $h(\mu)\geq\theta$ \ and \
$\vmid{A\cap[\mu,\mu+h(\mu))} = h(\mu).$ \ Now let \ $J_h$ \ be the
set of all \ $A\subseteq\kappa$ \ such that \ $T_\theta^A$ \ is a
nonstationary subset of \ $\kappa$ \ for some \ $\theta\in Y.$ \ It
is simple to check that \ $J_h$ \ is an ideal on \ $\kappa.$
\medskip
Let us remark in passing that if \ $\kappa$ \ is weakly Mahlo and \
$h$ \ is defined by~: $h(\mu) = \omega$ \ if \ $\mu$ \ is singular,
and \ $h(\mu) = \mu$ \ otherwise, then a subset \ $A$ \ of \
$\kappa$ \ lies in \ $J_h$ \ if and only if the set of all \ $\mu\in
Y$ \ such that \ $\mu$ \ is regular and \ $\vmid{A\cap[\mu,\mu+\mu)}
= \mu$ \ is nonstationary in \ $\kappa.$
\medskip
Let us show that \ $J_h$ \ is a semi-$Q$-point. Thus fix a \
$<\kappa$-to-one function \ $f : \kappa\hfltoutpetit{}{}\kappa.$ \ Then
$$C = \Bigr\{\mu\in Y : \mu =
\bigcup_{\alpha<\mu}f^{-1}(\{\alpha\})\Bigl\}$$
is a closed unbounded subset of \ $\kappa.$ \ Set \ $Q =
\displaystyle\bigcup_{\mu\in C}[\mu,\mu+h(\mu)).$ \ It is immediate that \
$\kappa-Q\in J_h.$ \ Now fix \ $\alpha\in\kappa$ \ such that \
$Q\cap f^{-1}(\{\alpha\})\not=\phi.$ \ Pick \ $\nu\in C$ \ so that
$$[\nu,\nu+h(\nu))\cap f^{-1}(\{\alpha\})\not=\phi.$$
Clearly, \ $\alpha\geq\nu$ \ and \ $\nu\cap f^{-1}(\{\alpha\}) =
\phi.$ \ Let \ $\rho$ \ be the least element of \ $C$ \ that is \
$>\nu.$ \ Then \ $\alpha<\rho$ \ and \
$f^{-1}(\{\alpha\})\subseteq\rho.$ \ Thus
$$Q\cap f^{-1}(\{\alpha\})\subseteq[\nu,\nu+h(\nu))$$
and consequently
$$\vmid{Q\cap f^{-1}(\{\alpha\})}\leq h(\nu)\leq\nu\leq\vmid\alpha.$$
\bigskip
It remains to show that \ $J$ \ is not weakly rapid. Fix \ $D\in
J_h^+.$ \ Then
$$S = \{\mu\in T_\omega^D : \vmid{T_\omega^D\cap\mu} = \mu\}$$
is a stationary subset of \ $\kappa.$ \ Given \ $\mu\in S, \
\vmid{D\cap\mu} = \mu$ \ since
$$D\cap[\rho,\rho+h(\rho))\subset D\cap\mu$$
for every \ $\rho\in\mu\cap T_\omega^D,$ \ and hence
$$\hbox{o.t.}(D\cap(\mu+h(\mu)) > \mu+1.$$
\remonte\hfill\carreblanc
\bigskip\bigskip
{\bf THEOREM 5.3.} \ {\sl Suppose \ $\kappa$ \ is a limit cardinal.
Then }\ ${\hbox{\nomgot U}}_\kappa\leq \hbox{\pmb{non}}_\kappa$(weak
semi-$Q$-point).
\bigskip
{\bf Proof.} \ Let \ $J$ \ be an ideal on \ $\kappa$ \ with \
$cof(J)<\hbox{\nomgot U}_\kappa.$ \ Let us show that \ $J$ \ is a
weak semi-$Q$-point. Thus fix \ $A\in J^+$ \ and a \
$<\kappa$-to-one function \ $f : A\hfltoutpetit{}{}\kappa.$ \
Select
$B_\beta\in J$ for $\beta<cof(J)$ so that $J =
\displaystyle\bigcup_{\beta<cof(J)}P(B).$ For $\beta<cof(J),$ define
$g_\beta\in{}^\kappa\kappa$ by~:
$$g_\beta(\alpha) = \hbox{ least element of }
(\bigcup_{\gamma>\alpha}f^{-1}(\{\gamma\}))-B_\beta.$$
There is $h\in{}^\kappa\kappa$ such that $\{\alpha\in\kappa : g_\beta(\alpha) =
h(\alpha)\} \not=\phi$ for every $\beta<cof(J).$ Define $C\subseteq ran(h)$
by~:
$h(\alpha)\in C$ just in case
$h(\alpha)\in\displaystyle\bigcup_{\gamma>\alpha}f^{-1}(\{\gamma\}).$ Then
clearly $C\in J^+\cap P(A).$ Moreover, $C\cap
f^{-1}(\{\alpha\})\subseteq\{h(\gamma) : \gamma<\alpha\}$ for every
$\alpha<\kappa.$\hfill$\carreblanc$
\bigskip\bigskip
{\bf THEOREM 5.4.} \ {\sl Suppose that \ $\kappa$ \ is a limit 
cardinal and \ $2^{<\kappa}
= \kappa.$ \ Then}
$$\overline{\hbox{\pmb{non}}}_\kappa\hbox{(weakly semiselective) } 
\leq \hbox{\nomgot
U}_\kappa\leq
\hbox{\pmb{non}}_\kappa \hbox{(weak }Q\hbox{-point).}$$
\bigskip
{\bf Proof.} \ The proof of the first inequality is an easy 
modification of that of Lemma 6.1 in
[MP1]
(which should be corrected by substituting ``$e\in[\omega]^{<\omega}$ such that
$B\subseteq\displaystyle\bigcup_{j\in e}\omega^{E_j}\cup\bigcup_{f\in z}B_f"$
for ``$e\in [\displaystyle\bigcup_{j\in\omega}\omega^{E_j}]^{<\omega}$ such
that \ $B\subseteq e\cup\displaystyle\bigcup_{f\in z}B_f"$). \ The 
second inequality is
proved as Proposition 5.3 in [MP1]. \hfill\carreblanc
\bigskip
{\bf Remark.} \ Suppose that \ $\kappa$ \ is a limit cardinal, \ 
$2^{<\kappa}=\kappa$ \
and \ $\hbox{\pmb{non}}_\kappa$(weakly semiselective) \ 
$<\kappa^{+\omega}.$ \ Then by
Proposition 4.1 and Theorems 2.2 and 5.4,
$$\hbox{\nomgot U}_\kappa = \hbox{\pmb{non}}_\kappa \hbox{(weakly 
selective) } =
\hbox{\pmb{non}}_\kappa
\hbox{(weakly semiselective).}$$

{\bf Remark.} \ It is consistent (see [MP1]) that \ $\hbox{\nomgot
U}_\omega<\hbox{\pmb{non}}_\omega$(weak \ $Q$-point), and that \
$\hbox{\pmb{non}}_\omega$(weak \
$Q$-point) $<\hbox{\pmb{non}}_\omega$(weak semi-$Q$-point). We do not 
know whether these
results can be generalized.
\bigskip
{\bf QUESTION.} \ Is it consistent that \ $\kappa$ \ is a limit cardinal, \
$2^{<\kappa}>\kappa$ \ and \ $\kappa^+<\hbox{\pmb{non}}_\kappa$(weak 
$Q$-point) ?
\bigskip
{\bf QUESTION.} \ By a result of [MP1], \ 
$cf(\hbox{\pmb{non}}_\omega$(weak $Q$-point))
$>\omega.$
\ Does this generalize ?

\vskip 1,5cm
{\Bbf{6. $\hbox{\pmb{non}}_\kappa (J^+\hfltoutpetit{}{} (J^+,\theta)^2)
$}}
\bigskip
In this section we use standard material to discuss the value of \ 
$\hbox{\pmb{non}}_\kappa
(J^+\hfltoutpetit{}{} (J^+,\theta)^2)$ \ for a cardinal \ $\theta\in 
[3,\kappa].$
\bigskip\bigskip
{\bf THEOREM 6.1.}
\medskip
{\parindent=1,5cm
\item{\bf (i)} $\hbox{\nomgot d}_\kappa \geq 
\hbox{\pmb{non}}_\kappa(J^+\hfltoutpetit{}{}
(J^+,3)^2).$
\medskip
\item{\bf (ii)} $\overline{\hbox{\nomgot 
d}}_\kappa\geq\overline{\hbox{\pmb{non}}}_\kappa
(J^+\hfltoutpetit{}{} (J^+,3)^2).$
\medskip
\item{\bf (iii)} $\overline{\hbox{\pmb{non}}}_\kappa$(weak $P$-point)
$\geq\overline{\hbox{\pmb{non}}}_\kappa (J^+\hfltoutpetit{}{}(J^+,\omega)^2).$
\par}
\bigskip
{\bf Proof.} \ (i) and (ii) : By a straightforward generalization of 
Lemma 4.4 in [M2],
there exists an ideal $J$ on $\kappa$ such that
$\overline{cof}(J)\leq\overline{\hbox{\nomgot d}}_\kappa, 
cof(J)\leq\hbox{\nomgot
d}_\kappa$ and
$J^+\hfltoutpetit{}{}\kern-3mm/\kern2mm(J^+,3)^2.$

(iii) : Baumgartner, Taylor and Wagon [BauTW] established that if \ $J$ \ is
an ideal on \ $\kappa$ \ such that \ $J^+\hfltoutpetit{}{} 
(J^+,\omega)^2,$ \ then \ $J$
\ is a weak \ $P$-point.~\hfill\carreblanc
\bigskip
{\bf Definition.} \ Given an ideal \ $J$ \ on \ $\kappa, \ A\in J^+$ 
\ and \ $F :
\kappa\times\kappa\hfltoutpetit{}{} 2, \ (J,A,F)$ \ is 0-good if 
there is \ $D\in J^+\cap
P(A)$ \ such that \ $\{\beta\in D : F(\alpha,\beta) = 1\}\in J$ \ for 
every \ $\alpha\in D.$
\bigskip
The following is readily checked.
\bigskip\bigskip
{\bf LEMMA 6.2.} \ {\sl Suppose that \ $J$ \ is weakly selective and 
\ $(J,A,F)$ \ is
0-good, where \ $J$ \ is an ideal  on \ $\kappa, \ A\in J^+$ \ and \ $F :
\kappa\times\kappa\hfltoutpetit{}{} 2.$ \ Then there is \ $B\in 
J^+\cap P(A)$ \ such that
\ $F$ \ is constantly 0 on \ $[B]^2.$}
\bigskip\bigskip
{\bf LEMMA 6.3.} \ {\sl Suppose that \ $(J,A,F)$ \ is not 0-good, 
where \ $J$ \ is an ideal
on \ $\kappa, \ A\in J^+$ \ and \ $F : 
\kappa\times\kappa\hfltoutpetit{}{}  2.$ \ Then~:}
\medskip
{\parindent=1,5cm
\item{\bf (i)} {\sl There is \ $B\subseteq A$ \ such that} \ o.t.$(B) 
= \omega+1$  \ {\sl
and
\
$F$ \ is identically 1 on \ $[B]^2.$}
\medskip
\item{\bf (ii)} {\sl Suppose that \ $\hbox{\nomgot a}_J>\kappa$ \ and 
\ $\theta$ \ is an
uncountable cardinal \ $<\kappa$ \ such that \
$\kappa\hfltoutpetit{}{}(\kappa,\theta)^2.$ \ Then there is \ 
$C\subseteq A$ \ such that}
\ o.t.$(C) = \theta+1$ \ {\sl and \ $F$ \ is identically 1 on \
$[C]^2.$}
\par}
\bigskip
{\bf Proof.} \ The proof is similar to that of Lemma 10.4 below. 
\hfill\carreblanc
\bigskip\bigskip
{\bf THEOREM 6.4.}
\medskip
{\parindent=1,5cm
\item{\bf (i)}
$\overline{\hbox{\pmb{non}}}_\kappa(J^+\hfltoutpetit{}{}(J^+,\omega+1)^2)\geq
\overline{\hbox{\pmb{non}}}_\kappa$(weakly selective).
\medskip
\item{\bf (ii)} {\sl Suppose that \ $\theta$ \ is an infinite 
cardinal \ $<\kappa$ \ such
that \ $\kappa\hfltoutpetit{}{}(\kappa,\theta)^2.$ \ Then}
$$\hbox{\bf{non}}_\kappa(J^+\hfltoutpetit{}{}(J^+,\theta+1)^2)\geq\hbox{\bf{non}}_\kappa
\hbox{(weakly selective).}$$
\par}
\bigskip
{\bf Proof.}
(i) : Baumgartner, Taylor and Wagon [BauTW] showed that \
$J^+\hfltoutpetit{}{}(J^+,\omega+1)^2$ \ for every weakly selective 
ideal \ $J$ \ on \
$\kappa.$
\medskip
(ii) : By Lemmas  6.2 and 6.3. \hfill\carreblanc
\bigskip
{\bf Remark.} \ Suppose that \ $\kappa$ \ is a successor cardinal and 
\ $\theta$ \ is
cardinal \ $\geq 2$ \ such that \ 
$\kappa\hfltoutpetit{}{}(\kappa,\theta)^2.$ \ Then by
Theorems 6.1 (i), 6.4 (ii) and 2.2 and Proposition 5.1, \ 
$\hbox{\nomgot{d}}_\kappa =
\hbox{\pmb{non}}_\kappa(J^+\hfltoutpetit{}{} (J^+,\theta+1)^2).$
\bigskip
{\bf Remark.} \ It is consistent (see [M2]) that \ $\hbox{\nomgot
d}>\hbox{\pmb{non}}_\omega(J^+\hfltoutpetit{}{}(J^+,3)^2).$ \ We do 
not know whether this
can be generalized.
\bigskip\bigskip
{\bf THEOREM 6.5.} \ {\sl Suppose \ $\kappa$ \ is a weakly compact 
cardinal. Then~:}
\medskip
{\parindent=1,5cm
\item{\bf (i)} $\overline{\hbox{\pmb{non}}}_\kappa$(weak\ $Q$-point)
$\geq\overline{\hbox{\pmb{non}}}_\kappa (J^+\hfltoutpetit{}{} (J^+,\kappa)^2).$
\medskip
\item{\bf (ii)} $\hbox{\pmb{non}}_\kappa(J^+\hfltoutpetit{}{}(J^+)^2) =
\hbox{\pmb{non}}_\kappa(J^+\hfltoutpetit{}{}(J^+,\kappa)^2) =
\hbox{\pmb{non}}_\kappa$(weakly selective).
\par}
\bigskip
{\bf Proof.} \ The result follows from Theorems 2.2 and 6.1 (i) and 
the following two
well-known facts : (1) Every ideal \ $J$ \ on \ $\kappa$ \ such that \
$J^+\hfltoutpetit{}{} (J^+,\kappa)^2$
\ is a weak \ $Q$-point ;  (2) If \ $\kappa$ \ is weakly compact, then \
$J^+\hfltoutpetit{}{} (J^+)^2$ \  for every weakly selective ideal \ 
$J$ \ on \ $\kappa$ \
such that\ $\hbox{\nomgot a}_J>\kappa.$
\hfill\carreblanc
\eject
\vskip 1,5cm
{\Bbf{7. $\hbox{\pmb{non}}_\kappa (J^+\hfltoutpetit{}{} [J^+]_\rho^2)$}}
\bigskip
In this section we consider the cardinal \ 
$\hbox{\pmb{non}}_\kappa(J^+\hfltoutpetit{}{}
[J^+]_\rho^2),$
\ where \ $3\leq\rho\leq\kappa,$ \ about which little is known. We 
begin with the case where
\ $\rho=3.$ \ The following is due to Blass [Bl].
\bigskip\bigskip
{\bf LEMMA 7.1.} \ {\sl Suppose \ $J$ \ is an ideal on \ $\kappa$ \ such that \
$J^+\hfltoutpetit{}{} [J^+]_3^2.$ \ Then \ $J$ \ is a weak \ $P$-point.
}
\bigskip
{\bf Proof.} \ Fix \ $A\in J^+$ \ and \ $f\in{}^A\kappa$ \ with \ 
$\{f^{-1}(\{\gamma\}) :
\gamma\in\kappa\}\subseteq J.$ \ Define \ $g : [A]^2\hfltoutpetit{}{} 
3$ \ by stipulating
that \
$g(\alpha,\beta) = 0$ \ if and only if \ $f(\alpha) < f(\beta),$ \ and \
$g(\alpha,\beta)=1$ \ if and only if \ $f(\alpha) = f(\beta).$ \ 
There are \ $B\in J^+\cap
P(A)$ \ and \ $i<3$ \ such that \ $i\notin g''[B]^2.$ \ It is simple 
to see that \
$i\not=0,$ \ so \ $f$ \ is \ $<\kappa$-to-one on \ $B.$ \hfill\carreblanc
\bigskip
The following is proved by adapting an argument of Baumgartner and 
Taylor [BauT].
\bigskip\bigskip
{\bf LEMMA 7.2.} \ {\sl Suppose \ $J$ \ is an ideal on \ $\kappa$ \ such that \
$J^+\hfltoutpetit{}{} [J^+]_3^2,$ \ and \ $(J,A,F)$ \ is 0-good, 
where \ $A\in J^+$ \ and \
$F : \kappa\times\kappa\hfltoutpetit{}{} 2.$ \ Then either there 
exists \ $C\in J^+\cap
P(A)$ \ such that \ $F$ \ is constantly 0 on \ $[C]^2,$ \ or for 
every \ $\delta<\kappa,$ \
there exists \ $Q\subseteq A$ \ such that  }\ o.t.$(Q) = \delta$ \ 
{\sl and \ $F$ \ is
constantly 1 on
\ $[Q]^2.$ }
\bigskip
{\bf Proof.} \ Select \ $B\in J^+\cap P(A)$ \ so that \ $\{\beta\in B 
: F(\alpha,\beta) =
1\}\in J$ \ for every \ $\alpha\in B.$ \ By Lemma 7.1, there exists \ 
$S\in J^+\cap P(B)$ \
so that \ $\vmid{\{\beta\in S : F(\alpha,\beta) = 1\}}<\kappa$ \ for 
every \ $\alpha\in S.$
\ Define \ $\delta_\xi$ \ for \ $\xi<\kappa$ \ by~:
\medskip
(i) $\delta_0 = \cap S ;$
\medskip
(ii) $\delta_{\xi+1} =$ the least \ $\zeta<\kappa$ \ with the property that \
$\zeta>\beta$ \ for every \ $\beta\in S$ \ such that \ 
$F(\alpha,\beta) = 1$ \ for some \
$\alpha\in S\cap\delta_\zeta$ ;
\medskip
(iii) $\delta_\xi = \displaystyle\bigcup_{\zeta<\xi}\delta_\zeta$ \ if \
$\xi$ \ is a limit ordinal \ $>0.$

\medskip
Let \ $X$ \ be the set of all limit ordinals \ $<\kappa.$ \ For \ 
$\eta\in X, \ n\in\omega$
\ and \ $j<2,$ \ set
$$d_{\eta,n}^j = S\cap [\delta_{\eta+2n+j},\delta_{\eta+2n+j+1}).$$
For \ $j<2,$ \ let
$$D^j = \cup\{d_{\eta,n}^j : \eta\in X \hbox{ and } n\in\omega\}.$$
Select \ $k<2$ \ so that \ $D^k\in J^+.$ \ Notice that \ 
$F(\alpha,\beta) = 0$ \ if \
$(\alpha,\beta)\in [D^k]^2$ \ and \ $\{\alpha,\beta\}\not\subseteq 
d_{\eta,n}^k$ \ for all
\ $\eta\in X$ \ and \ $n\in\omega.$
\bigskip
Define \ $h : [D^k]^2 \hfltoutpetit{}{} 3$ \ by stipulating that \ 
$h(\alpha,\beta) = 0$ \
if and only if \ $\{\alpha,\beta\}\not\subseteq d_{\eta,n}^k$ \ for 
all \ $\eta\in X$ \ and
\
$n\in\omega,$ \ and \ $h(\alpha,\beta) = 1$ \ if and only if \ 
$F(\alpha,\beta) = 1.$ \
There are \ $W\in J^+\cap P(D^k)$ \ and \ $i<3$ \ so that \ $i\notin 
h''[W]^2.$ \ Clearly,
\ $i\not=0.$ \ If \ $i=1, \ F$ \ is identically 0 on \ $[W]^2.$ \ Now 
assume  \ $i=2.$ \
Let \ $Z$ \ be the set of all \ $(\eta,n)\in X\times\omega$ \ such that \
$W\cap d_{\eta,n}^k\not=\phi.$ \ Suppose that there is \ 
$\gamma<\kappa$ \ such that \
o.t.$(W\cap d_{\eta,n}^k)\leq\gamma$ \ for every \ $(\eta,n)\in Z.$ \ 
Then there exists \
$C\in J^+\cap P(W)$ \ such that \ $\vmid{C\cap d_{\eta,n}} = 1$ \ for 
any \ $(\eta,n)\in
Z.$ \ Clearly, \ $F$ \ takes the constant value 0 on \ $[T]^2.$ 
\hfill\carreblanc
\bigskip\bigskip
{\bf PROPOSITION 7.3.} \ {\sl Suppose \ $\theta\in(2,\kappa)$ \ {\sl 
is a cardinal such
that } \ $\kappa\hfltoutpetit{}{}(\kappa,\theta)^2.$ \ Then }\
$$\hbox{\pmb{non}}_\kappa(J^+\hfltoutpetit{}{} [J^+]_3^2)\leq
\hbox{\pmb{non}}_\kappa(J^+\hfltoutpetit{}{} (J^+,\theta+1)^2).$$
\bigskip
{\bf Proof.} \ By Theorem 2.2 and Lemmas 6.3, 7.1 and 7.2. \hfill\carreblanc
\bigskip
Let us now consider the partition relation \ $J^+\hfltoutpetit{}{} 
[J^+]_\kappa^2.$ \ We
begin with the following observation.
\bigskip\bigskip
{\bf PROPOSITION 7.4.} \ {\sl Suppose \ $\kappa$ \ is inaccessible. 
Then there is an ideal
\ $J$ \ on \ $\kappa$ \ such that} (a) \ 
$J^+\hfltoutpetit{}{}\kern-3mm/\kern2mm
[J^+]_\kappa^2,$ \ (b) $J$ \ {\sl is not a weak semi}-$Q$-point,  (c) 
$\hbox{\nomgot
a}_J>\kappa,$ {\sl and} (d) $J^+\hfltoutpetit{}{} (J^+,\alpha)^2$ \ for every \
$\alpha<\kappa.$
\bigskip
{\bf Proof.} \ Let \ $<\rho_\alpha : \alpha<\kappa>$ \ be the 
increasing enumeration of all
strong limit infinite cardinals \ $<\kappa.$ \ let \ $Z$ \ be the set 
of all regular
infinite cardinals \ $<\kappa.$ \ For \ $\mu\in Z,$ \ set \ $\nu_\mu 
= (\rho_\mu)^{++}.$ \
Then \ 
$\nu_\mu\hfltoutpetit{}{}\kern-3mm/\kern2mm[\nu_\mu]_{\nu_\mu}^2$ \ 
by a result of
Todorcevic [To1]. On the other hand, by a result of Erd\"os and Rado 
(see [EHM\'aR], Corollary
17.5), \ $\nu_\mu\hfltoutpetit{}{} (\nu_\mu,\tau)^2$ \ for every 
infinite cardinal \
$\tau<\mu.$ \ Pick pairwise disjoint \ $A_\mu$ \ for \ $\mu\in Z$ \ so that \
$\vmid{A_\mu}=\nu_\mu$ \ for any \ $\mu\in Z,$ \ and \ 
$\displaystyle\bigcup_{\mu\in Z}A_\mu
=\kappa.$ \ Let  \ $J$ \ be the set of all \ $B\subseteq\kappa$ \ such that
$$\vmid{\{\mu\in Z : \vmid{B\cap A_\mu} = \nu_\mu\}}<\kappa.$$
It is simple to see that \ $J$ \ is an ideal on \ $\kappa.$
\bigskip
For \ $\mu\in Z,$ \ pick \ $g_\mu : 
[A_\mu]^2\hfltoutpetit{}{}\nu_\mu$ \ so that \
$g''_\mu[B]^2 = \nu_\mu$ \ for every \ $B\subseteq A_\mu$ \ with \ 
$\vmid B=\nu_\mu.$ \ Let
\ $G : [\kappa]^2\hfltoutpetit{}{}\kappa$ \ be such that \ 
$\displaystyle\bigcup_{\mu\in
Z}g_\mu\subseteq G.$ \ Then clearly \ $G''[C]^2 = \kappa$ \ for any \ 
$C\in J^+.$

Define \ $f\in{}^\kappa\kappa$ \ by stipulating that \ $f^{-1}(\{\mu\}) 
= A_\mu$ \ for
every \ $\mu\in Z.$ \ Clearly, there is no \ $S\in J^+$ \ so that \ 
$\vmid{S\cap
f^{-1}(\{\alpha\}}\leq\vmid\alpha$ \ for all \ $\alpha<\kappa.$ \ 
Hence \ $J$ \ is not a
weak semi-$Q$-point.
\bigskip
Let us next show that \ $\hbox{\nomgot a}_J>\kappa.$ \ Thus suppose 
that \ $B_\alpha\in
J^+$ \ for \ $\alpha<\kappa,$ \ and \ $B_\alpha\cap B_\beta\in J$ \ whenever  \
$\beta<\alpha<\kappa.$ \ Select a strictly increasing function \ $k :
\kappa\hfltoutpetit{}{} Z$ \ so that
$$\vmid{(B_\alpha - (\bigcup_{\beta<\alpha}B_\beta))\cap 
A_{k(\alpha)}} = \nu_{k(\alpha)}$$
for any \ $\alpha<\kappa.$ \ Set
$$T = \bigcup_{\alpha<\kappa}((B_\alpha-(\bigcup_{\beta<\alpha}B_\beta))\cap
A_{k(\alpha)}).$$
Then \ $T\in J^+$ \ and moreover  \ $\vmid{T\cap B_\alpha}<\kappa$ \ 
for every \
$\alpha<\kappa.$
\bigskip
It remains to prove (d). Thus fix \ $A\in J^+$ \ and \ $F :
\kappa\times\kappa\hfltoutpetit{}{} 2.$ \ Suppose that
there is  \ $\eta<\kappa$ \ such that for every \ $Q\subseteq A$ \ 
with \ o.t.$(Q) = \eta,
\ F$ is not constantly 1 on \ $[Q]^2.$ \ Since by Theorem 17.1 of 
[EHM\'aR] \ $\kappa\hfltoutpetit{}{} (\kappa,\alpha)^2$ \ for every
\ $\alpha<\kappa,$ \ it follows from Lemma 6.3 that \ $(J,A,F)$ is $0$-good.
  Select \ $D\in J^+\cap P(A)$ \ so that \
$\{\beta\in D : F(\alpha,\beta) = 1\}\in J$ \ for every \ $\alpha\in 
D.$ \ Define \
$D_\gamma$ \ for \
$\gamma<\kappa$ \ and a strictly increasing function \ $h : 
\kappa\hfltoutpetit{}{} Z$ \ so
that
\medskip
(0) \hskip 1,5cm $D_\gamma = D-(\displaystyle\bigcup_{\delta<\gamma}\ 
\bigcup_{\alpha\in
D_\delta\cap A_{h(\delta
)}}\{\beta\in D : F(\alpha,\beta) = 1\})$ ;

(1) \hskip 1,5cm $\vmid{D_\gamma\cap A_{h(\gamma)}} = \nu_{h(\gamma)}.$

\medskip
For \ $\gamma\in(\vmid\eta^+,\kappa),$ \ select \ $X_\gamma\subseteq 
D_\gamma\cap
A_{h(\gamma)}$ \ so that \ $\vmid{X_\gamma} = \nu_{h(\gamma)}$ \ and 
\ $F$ \ is constantly
0 on \ $[X_\gamma]^2.$ \ Set \ $Y = \displaystyle\bigcup_{\mid
\eta\mid^+<\gamma<\kappa}X_\gamma.$ \ Then clearly \ $Y\in J^+\cap 
P(A).$ \ Moreover, \ $F$
\ takes the constant value 0 on \
$[Y]^2.$ \hfill \carreblanc
\bigskip

{\bf Remark.} \ $J^+\hfltoutpetit{}{} (J^+,\kappa)^2$ \ does not 
necessarily imply that \
$J^+\hfltoutpetit{}{}  [J^+]_\kappa^2.$ \ This follows from the 
following two facts~: (0)
If \
$\kappa$ \ is weakly compact, then there exists a normal ideal \ $J$ 
\ on \ $\kappa$ \ such
that \ $J^+\hfltoutpetit{}{} (J^+,\kappa)^2$ \ ([Bau1], [Bau2]) ; (1) 
Assuming \ $V=L, \
\kappa$ \ is completely ineffable if and only if there is a normal 
ideal \ $J$ \ on \
$\kappa$ \ such that \ $J^+\hfltoutpetit{}{} [J^+]_\kappa^2$ \ ([M4]).
\bigskip
Recall that for  \ $S\subseteq\kappa, \lozenge_\kappa^*(S)$ \ means 
that there are \
$s_\alpha\in P_{{\mid\alpha\mid^+}}(\alpha)$ \ for \ $\alpha\in S$ \ 
such that for every \
$A\subseteq\kappa,$ \ there exists a closed unbounded subset \ $C$ \ 
of \ $\kappa$ \ with
the property that \ $A\cap\alpha\in s_\alpha$ \ for every \ 
$\alpha\in C\cap S.$
\bigskip\bigskip
{\bf PROPOSITION 7.5.-} \ {\sl Suppose that \ $\lozenge_\kappa^*(S)$ 
\ holds for some
stationary subset \ $S$ \ of \ $\kappa.$ \ Then }\ $\hbox{\nomgot
d}_\kappa\geq\hbox{\pmb{non}}_\kappa(J^+\hfltoutpetit{}{} 
[J^+]_\kappa^2)$ \ {\sl and } \
$\overline{\hbox{\nomgot
d}}_\kappa\geq\overline{\hbox{\pmb{non}}}_\kappa(J^+\hfltoutpetit{}{} 
[J^+]_\kappa^2).$
\bigskip
{\bf Proof.} \ By a result of [M4], the hypothesis implies that \
$NS_\kappa^+\hfltoutpetit{}{}\kern-3mm/\kern2mm 
[NS_\kappa^+]_\kappa^2.$ \hfill\carreblanc\eject
\bigskip
{\bf Remark.} \ It is shown in [S] that if (a) \ $\kappa$ \ is a 
successor cardinal \
$\geq\omega_2$ \ with \ $2^{<\kappa}=\kappa,$ \ and (b) setting \ 
$\kappa=\nu^+, \
\mu^\tau\leq\nu$ \ for every infinite cardinal \ $\mu<\nu,$ \ where \ 
$\tau = \aleph_1$ \
if \ $cf(\nu)=\omega$ \ and \ $\tau=\aleph_0$ \ otherwise, then there 
is a stationary
subset \ $S$ \ of \ $\kappa$ \ such that \ $\lozenge_\kappa^*(S)$ \ holds.
\bigskip
{\bf Remark.} \ We do not know whether it is consistent that the 
conclusion of Proposition
7.5 fails. Results of Section 15 (below) imply that
$$\hbox{\pmb{non}}_\kappa(J^+\hfltoutpetit{}{} [J^+]_\kappa^2)\leq 
(\overline{\hbox{\nomgot
d}}_\kappa)^{<\kappa}$$
if \ $\kappa$ \ is a limit cardinal such that \ $2^{<\kappa}=\kappa.$
\vskip 1,5cm
{\Bbf{8. $\hbox{\nomgot d}_{\kappa,\lambda}^\kappa$}}
\bigskip
We now start our study of combinatorial properties of ideals on \ 
$P_\kappa(\lambda).$ \
The aim of this section is to present a two-cardinal version of Theorem 2.2.
\bigskip
\bigskip
{\bf Definition.} \ $\hbox{\nomgot d}_{\kappa,\lambda}^\kappa$ \ is 
the least cardinality
of any \ $F\subseteq{}^\kappa(P_\kappa(\lambda))$ \ with the property 
that for every \break
$g\in{}^\kappa(P_\kappa(\lambda)),$ \ there is \ $f\in F$ \ such that 
\ $g(\alpha)\subseteq
f(\alpha)$ \ for all \ $\alpha\in\kappa.$
\bigskip
{\bf Remark.} \ It is shown in [MP\'eS1] that \ $\hbox{\nomgot 
d}_{\kappa,\lambda}^\kappa =
\hbox{\nomgot d}_\kappa\cdot u(\kappa^+,\lambda).$
\bigskip
{\bf Definition.} \ Given an ideal \ $H$ \ on \ $P_\kappa(\lambda), \ {\cal
M}_H^{\geq\kappa}$\ is the set of all \ $Q\subseteq H^+$ \ such that 
(i) \ $\vmid
Q\geq\kappa,$ \ (ii) $A\cap B\in H$ \ for all \ $A,B\in Q$ \ with \ 
$A\not= B,$ \ and (iii)
for every \ $C\in H^+,$ \ there is \ $A\in Q$ \ with \ $A\cap C\in H^+.$

$\hbox{\nomgot a}_H$ \ is the least cardinality of any member of \ ${\cal
M}_H^{\geq\kappa}$ \ if \ ${\cal M}_H^{\geq\kappa}\not=\phi,$ \ and \
$2^{(\lambda^{<\kappa})^+}$ \ otherwise.
\bigskip
The following is proved as Proposition 11.2 of [MP2].

\bigskip\bigskip
{\bf PROPOSITION 8.1.}\ {\sl Given a $\kappa$-normal ideal $H$ on
$P_\kappa(\lambda),$ the following are equivalent~:
\medskip
{\parindent=1,5cm
\item{\bf (i)} $\hbox{\nomgot a}_H = \kappa.$
\item{\bf (ii)} $sat(H) >\kappa.$
\par}}
\bigskip\bigskip
{\bf COROLLARY 8.2.}\ {\sl Let $A\in(NS_{\kappa,\lambda}^\kappa)^+$ and
set $H = NS_{\kappa,\lambda}^\kappa\mid A.$ Then $\hbox{\nomgot a}_H
= \kappa.$}
\bigskip
{\bf Proof.}\ The result follows from Proposition 8.1 since
$sat(H)>\kappa$ by a result of Abe [A]. \hfill$\carreblanc$
\bigskip
The following is proved as Proposition 11.1 (ii) of [MP2].
\bigskip\bigskip
{\bf PROPOSITION 8.3.}\ {\sl Given an ideal $H$ on $P_\kappa(\lambda),$
the following are equivalent~:
\medskip
{\parindent=1,5cm
\item{\bf (i)} $\hbox{\nomgot a}_H = \kappa.$
\item{\bf (ii)} There exist $A_\alpha\in H^+$ for $\alpha<\kappa$
such that (a) $A_\alpha\subseteq A_\beta$ whenever
$\beta<\alpha<\kappa,$ and (b) for every $C\in H^+,$ there is
$\alpha<\kappa$ such that $C-A_\alpha\in H^+.$
\par}}
\bigskip\bigskip
{\bf Definition.} \ An ideal \ $H$ \ on \ $P_\kappa(\lambda)$ \ is a 
weak \ $\pi$-point if
given \ $f\in {}^\kappa H$ \ and \ $A\in H^+,$ \ there is \ $B\in 
H^+\cap P(A)$ \ such that
\ $B\cap f(\alpha)\in I_{\kappa,\lambda}$ \ for every \ $\alpha\in\kappa.$
\bigskip\bigskip
{\bf THEOREM 8.4.} \ {\sl Let \ $H$ \ be an ideal on \ 
$P_\kappa(\lambda)$ \ such that \
$cof(H) <\hbox{\nomgot d}_{\kappa,\lambda}^\kappa.$ \ Then \ 
$\hbox{\nomgot a}_H>\kappa$ \
and \ $H$ \ is a weak \ $\pi$-point.
}
\bigskip
{\bf Proof.} \ Let $A_\alpha\in H^+$ for $\alpha<\kappa$ be such that
$A_\alpha\subseteq A_\beta$ for all $\beta<\alpha.$ Select $X\subseteq
H$ so that $\vmid X = cof(H)$ and $H = \displaystyle\bigcup_{B\in
X}P(B).$ For $B\in X,$ define $f_B\in{}^\kappa(P_\kappa(\lambda))$ so
that $f_B(\alpha)\in A_\alpha-B.$ There is
$g\in{}^\kappa(P_\kappa(\lambda))$ such that $\{\alpha<\kappa :
g(\alpha)\not\subseteq f_B(\alpha)\}\not=\phi$ for every $B\in X.$ Set
\ $C = \displaystyle\bigcup_{\alpha\in\kappa}\{a\in A_\alpha :
g(\alpha)\not\subseteq a\}.$ \
Then $C\in H^+,$ and moreover $C-A_\alpha\in I_{\kappa,\lambda}$ for
any $\alpha<\kappa.$

\bigskip
{\bf Definition.} \ $\overline{\hbox{\nomgot
d}}_{\kappa,\lambda}^\kappa$ \ is the least cardinality of any \
$X\subseteq {}^\kappa(P_\kappa(\lambda))$ \  with the property \
that for every \ $g\in {}^\kappa(P_\kappa(\lambda)),$ \ there is
$x\in P_\kappa(X)$ \ such that \
$g(\alpha)\subseteq\displaystyle\bigcup_{f\in x} f(\alpha)$
for every \ $\alpha<\kappa.$
\bigskip
{\bf Remark.} \ It is shown in [MRoS] that \ ${\overline{\hbox{\nomgot
d}}_{\kappa,\lambda}^\kappa} = \overline{\hbox{\nomgot
d}}_\kappa\cdot\hbox{cov}(\lambda,\kappa^+,\kappa^+,\kappa),$ \ where \
$\hbox{cov}(\lambda,\kappa^+,\kappa^+,\kappa)$ \ denotes the least
cardinality of any \ $X\subseteq P_{\kappa^+}(\lambda)$ \ such that
for every \ $b\in P_{\kappa^+}(\lambda),$ \ there is \ $x\in
P_\kappa(X)$ \ with \ $b\subseteq\cup x.$
\bigskip
{\bf Remark.} \ It is immediate that \ $I_{\kappa,\lambda}$ \ is a 
weak \ $\pi$-point. On
the other hand \ $\hbox{\nomgot a}_{I_{\kappa,\lambda}}>\kappa$ \ 
does not necessarily
hold. In fact if \ $cf(\lambda)\not=\kappa$ \ and \ $\overline{\hbox{\nomgot
d}}_{\kappa,\sigma}^\kappa\leq\lambda$ \ for every cardinal \ 
$\sigma\in[\kappa,\lambda),$ \
then \ $\hbox{\nomgot a}_{I_{\kappa,\lambda}}=\kappa$ \ ([M6]).
\bigskip\bigskip
{\bf LEMMA 8.5.} \ {\sl Suppose that \ $H$ \ is an ideal on \ 
$P_\kappa(\lambda)$ \ with \
$\hbox{\nomgot a}_H = \kappa.$ \ Then there is an ideal \ $K$ \ on \ 
$P_\kappa(\lambda)$ \
such that (a) \ $K$ \ is not a weak \ $\pi$-point, (b) $cof(K)\leq 
cof(H),$ and (c)
$\overline{cof}(K) \leq\overline{cof}(H).$
}
\bigskip
{\bf Proof.} \ Select \ $A_\alpha\in H^+$ \ for \ $\alpha<\kappa$ \ 
so that ($\alpha$)
$A_\alpha\subseteq A_\beta$ \ whenever \ $\beta<\alpha<\kappa,$ \ and 
($\beta$) for any \
$C\in H^+,$ \ there is \ $\alpha<\kappa$ \ with \ $C-A_\alpha\in 
H^+.$ \ Let \ $K$ \ be the
set of all \ $B\subseteq P_\kappa(\lambda)$ \ such that \ $B\cap 
A_\alpha\in H$ \ for some
\ $\alpha<\kappa.$ \ It is simple to check that \ $K$ \ is as 
desired. \hfill\carreblanc
\bigskip\bigskip
{\bf THEOREM 8.6.}
\medskip
{\parindent=1,5cm
\item{\bf (i)} {\sl There is an ideal \ $H$ \ on \ 
$P_\kappa(\lambda)$ \ such that (a)
$\hbox{\nomgot a}_H=\kappa,$ \ (b) $cof(H) = 
\hbox{\nomgot{d}}_{\kappa,\lambda}^\kappa,$ \
  and }(c) \
$\overline{cof}(H)
\leq
\overline{\hbox{\nomgot d}}_{\kappa,\lambda}^\kappa.$
\medskip
\item{\bf (ii)} {\sl There is an ideal \ $K$ \ on \ 
$P_\kappa(\lambda)$ \ such that (a) $K$
\ is not a weak \ $\pi$-point, (b) $cof(K) = \hbox{\nomgot 
d}_{\kappa,\lambda}^\kappa,$ \
and (c) $\overline{cof}(K) \leq\overline{\hbox{\nomgot 
d}}_{\kappa,\lambda}^\kappa.$}
\par}
\bigskip
{\bf Proof.} \
(i) : Set \ $H = NS_{\kappa,\lambda}^\kappa.$ \ Then \ $\hbox{\nomgot a}_H =
\kappa$ \ by Corollary 8.2. Moreover,  \ $cof(H) = \hbox{\nomgot 
d}_{\kappa,\lambda}^\kappa$
([MP\'eS1]) and \ $\overline{cof}(H) = \overline{\hbox{\nomgot 
d}}_{\kappa,\lambda}^\kappa$
([MRoS]).
\medskip
(ii) : By (i), Lemma 8.5 and Theorem 8.4. \hfill\carreblanc

\bigskip
{\bf Remark.} \ Theorem 8.6 is not optimal, even under GCH. In
fact, suppose that the GCH holds, $\lambda = \sigma^+,$ where
$\sigma$ is a cardinal of cofinality $<\kappa,$ and $\kappa$ is not
the successor of a cardinal of cofinality $\leq cf(\sigma).$ Then \
$\overline{\hbox{\nomgot d}}_{\kappa,\lambda}^\kappa=\lambda$ 
([MRoS]). Moreover, there is\
$A\in(NS_{\kappa,\lambda}^\kappa)^+$ such that
$\overline{cof}(NS_{\kappa,\lambda}^\kappa\mid A) = \sigma$ ([MP\'eS2]).
Hence there is by Corollary 8.2 an ideal $H$ on
$P_\kappa(\lambda)$ (namely $H = NS_{\kappa,\lambda}^\kappa\mid A)$
such that $\overline{cof}(H) < \overline{\hbox{\nomgot
d}}_{\kappa,\lambda}^\kappa$ and $\hbox{\nomgot a}_H = \kappa,$ \ and 
by Lemma 8.5 an ideal
\ $K$ \ on \ $P_\kappa(\lambda)$ \ such that \ $\overline{cof}(K) 
<\overline{\hbox{\nomgot
d}}_{\kappa,\lambda}^\kappa$ \ and \ $K$ \ is not a weak \ $\pi$-point.
\vskip 1,5cm
{\Bbf{9. Weak \ $\chi$-pointness}}
\bigskip
{\bf Definition.} An ideal $H$ on $P_\kappa(\lambda)$ is a weak
$\chi$-point if given $A\in H^+$ and
$g\in{}^\kappa(P_\kappa(\lambda)),$ there is \break $B\in H^+\cap P(A)$ such
that $g(\cup(a\cap\kappa)\subseteq b$ for all $a,b\in B$ with
$\cup(a\cap\kappa)<\cup(b\cap\kappa).$
\bigskip
Our primary concern in this section is with the problem of determining when \
$I_{\kappa,\lambda}$ \ is a weak \ $\chi$-point. \break We will first 
give a sufficient
condition and then prove that this condition is  necessary if \ 
$\kappa$ \ is inaccessible.
\bigskip
The following is proved as Lemma 2.1 in [M2].
\bigskip\bigskip
{\bf THEOREM 9.1.}\ {\sl Let $H$ be an ideal on $P_\kappa(\lambda)$ such that}
$cof(H) < \hbox{{\bf cov}}(\pmb{M}_{\kappa,\kappa}).$ {\sl Then $H$ is a weak
$\chi$-point.}
\bigskip\bigskip
{\bf QUESTION.}\ Is it consistent that $2^{<\kappa}>\kappa$ and
$I_{\kappa,\kappa^+}$ is a weak $\chi$-point ?
\bigskip\bigskip
{\bf THEOREM 9.2.}\ {\sl Suppose that for all \ $A\in 
I_{\kappa,\lambda}^+$ \ with \
$A\subseteq\{a:\cup(a\cap\kappa)\in a\},$
\ there is \ $B\in
I_{\kappa,\lambda}^+\cap P(A)$ \ such that \
$\cup(a\cap\kappa)\in b$ \ for all \ $a,b\in B$ \ with \
$\cup(a\cap\kappa)<\cup(b\cap\kappa).$ \ Then} \
$\sigma<\overline{\hbox{\pmb{non}}}_\kappa$(weakly selective) {\sl for every \
$\sigma\in{\cal K}(\kappa,\lambda).$ }
\bigskip
{\bf Proof.} \ Suppose that \ $T\subseteq P_\kappa(\lambda-\kappa)$ \ 
is such that \
$\vmid{T\cap P(a)}<\kappa$ \ for every \ $a\in P_\kappa(\lambda),$ \ 
and \ $J$ \ is an
ideal on \ $\kappa$ \ with \ $\overline{cof}(J)\leq\vmid T.$ \ Select 
\ $D_d\in J$ \
for \ $d\in T$ \ so that for every \ $W\in J,$ \ there is \ $u\in 
P_\kappa(T)-\{\phi\}$ \
with \ $W\subseteq\displaystyle\bigcup_{d\in u}D_d.$ \  Now fix \ 
$G_\alpha\in J$ \ for \
$\alpha<\kappa.$ \ Define \ $A\subseteq P_\kappa(\lambda)$ \ by 
stipulating that \ $a\in A$ \
if and only if there is \
$\delta<\kappa$ \ such that (a) $\delta=\hbox{max}(a\cap\kappa),$ (b) \
$\delta\notin\displaystyle\bigcup_{d\in T\cap  P(a)}D_d,$ \ and (c) 
$\delta\notin G_\alpha$
\ for every \ $\alpha\in a\cap\delta.$
\bigskip
Let us show that \ $A\in I_{\kappa,\lambda}^+.$ \ Given \ $c\in 
P_\kappa(\lambda),$ \ pick
\ $\delta<\kappa$ \ so that \
$\delta\notin\displaystyle\bigcup_{d\in T\cap P(c)}D_d$ \ and for 
every \ $\alpha\in c\cap\kappa,
\delta>\alpha$ \ and \ $\delta\notin G_\alpha.$ \ Set \ 
$e=c\cup\{\delta\}.$ \ Then \ $e\in A.$
\bigskip
By our assumption there is \ $B\in I_{\kappa,\lambda}^+\cap P(A)$ \ such that \
$\cup(a\cap\kappa)\in b$ \ for all \ $a,b\in B$ \ with \
$\cup(a\cap\kappa)<\cup(b\cap\kappa).$ \ Set \ $C=\{\cup(a\cap\kappa) 
: a\in B\}.$ \ Then \ $C\in
J^+.$ \ Moreover, \ $\xi\notin G_\zeta$ \ for all \ $\zeta,\xi\in C$ 
\ with \ $\zeta<\xi.$
\hfill\carreblanc

\bigskip
We mention the following partial converse to Theorem 9.2.
\bigskip\bigskip
{\bf PROPOSITION 9.3.} \ {\sl Suppose  that \ $2^{<\kappa}=\kappa$ \ 
and \ $H$ \ is an
ideal on \ $P_\kappa(\lambda)$ \ such that }\ $cof(H) < 
\hbox{\pmb{non}}_\kappa$(weakly
selective). {\sl Then for all \ $f\in {}^\kappa\kappa$ \ and \ $A\in 
H^+,$ \ there is \
$B\in H^+\cap P(A)$ \ such that \ $f(\cup(a\cap\kappa))\subseteq b$ \ 
for all \ $a,b\in B$
\ with
\ $\cup(a\cap\kappa)<\cup(b\cap\kappa).$
  }
\bigskip
{\bf Proof.} \ Fix \ $f\in{}^\kappa\kappa$ \ and \ $A\in H^+.$ \ For 
\ $D\subseteq
P_\kappa(\kappa),$ \ set \ $Z_D = \{a\in P_\kappa(\lambda) : 
a\cap\kappa\in D\}.$ \ It is
simple to see that (a) $Z_{P_\kappa(\kappa)}=P_\kappa(\lambda),$ (b) 
$Z_{\cup{\hbox{\nomgot
D}}}Z_D$ \ for \ $\hbox{\nomgot D}\subseteq P(P_\kappa(\kappa)),$ (c) $Z_D\in
I_{\kappa,\lambda}$ \ for every \ $D\subseteq P_\kappa(\kappa)$ \ 
with \ $\vmid D = 1,$ \
and (d) $Z_{D'}\subseteq Z_D$ \ for all \ $D, D'\subseteq 
P_\kappa(\kappa)$ \ such that \
$D'\subseteq D.$ \ Hence
$$K = \{D\subseteq P_\kappa(\kappa) : Z_D\in H\mid A\}$$
is a \ $\kappa$-complete ideal on \ $P_\kappa(\kappa).$ \ For \ $C\subseteq
P_\kappa(\lambda),$ let \ $W_C$ \ be the set of all \ $d\in 
P_\kappa(\kappa)$ \ such that
$$\{a\in P_\kappa(\lambda) : a\cap \kappa = d\}\subseteq C.$$
If \ $C\in H\mid A,$ \ then \ $W_C\in K$ \ since \ $Z_{W_C}\subseteq 
C.$ \ Moreover, if \
$D\subseteq P_\kappa(\kappa)$ \ and \ $Z_D\subseteq C\subseteq 
P_\kappa(\lambda),$ \ then \
$D\subseteq W_C.$ \ Hence
$$cof(K) \leq cof(H\mid A)\leq cof(H).$$
For \ $d\in P_\kappa(\kappa),$ \ let \ $S_d$ \ be the set of all \ 
$e\in P_\kappa(\kappa)$
\ such that \ $f(\cup d)\not\subseteq e$ \ or \ $\cup e\leq\cup d.$ \ 
Then \ $S_d\in K$ \
since
$$\{a\in Z_{S_d} : f(\cup d)\cup\{(\cup d)+1\}\subseteq a\}=\phi.$$
Select a bijection \ $\ell : 
P_\kappa(\kappa)\hfltoutpetit{}{}\kappa.$ \ Since \
$cof(K)<\hbox{\bf non}_\kappa$ (weakly selective), there is \ $D\in 
K^+$ \ such that \
$e\notin S_d$ \ for all \ $d,e\in D$ \ such that \ $\ell(d)<\ell(e).$ \ Set
$$B =A\cap Z_D = \{a\in A : a\cap \kappa\in D\}.$$
Then \ $B\in H^+.$ \ Now fix \ $a,b\in B$ \ with \ 
$\cup(a\cap\kappa)<\cup(b\cap\kappa).$
\ Then clearly \ $\ell(a\cap\kappa)\not=\ell(b\cap\kappa).$ \ In fact \
$\ell(a\cap\kappa)<\ell(b\cap\kappa)$ \ (since otherwise \
$a\cap\kappa\notin S_{b\cap\kappa}$ \ and therefore \
$\cup(a\cap\kappa)>\cup(b\cap\kappa)).$ \ Hence \ $b\cap\kappa\notin 
S_{a\cap\kappa},$ \ so
\ $f(\cup(a\cap\kappa))\subseteq b\cap\kappa.$\hfill\carreblanc
\bigskip
{\bf Definition.} \ For \ $A\subseteq P_\kappa(\lambda),$ \ let
$$[A]_\kappa^2 = \{(\cup(a\cap\kappa),b) : a,b\in A \hbox{ and }
\cup(a\cap\kappa) < \cup(b\cap\kappa)\}$$
{\bf Remark.}
$$[P_\kappa(\lambda)]_\kappa^2 =
\{(\alpha,b)\in\kappa\times P_\kappa(\lambda) : \alpha < \cup(b\cap\kappa)\}.$$
\bigskip
{\bf Definition.} \ For \  $a,b\in P_\kappa(\lambda),$ \ let \ 
$a\prec b$ \ just in case \
$a\subseteq b$ \ and \ $\cup(a\cap\kappa) <\cup(b\cap\kappa).$
\bigskip\bigskip

{\bf Definition.} \ For \ $A\subseteq P_\kappa(\lambda),$ \ let

$$[A]_\prec^2 = \{(\cup(a\cap\kappa),b) : a,b\in A \hbox{ and } a\prec b\}.$$
\bigskip
{\bf Remark.} \ $[P_\kappa(\lambda)]_\prec^2 = [P_\kappa(\lambda)]_\kappa^2.$

\bigskip\bigskip
{\bf THEOREM 9.4.}\ {\sl Suppose that \ $\kappa$ \ is inaccessible 
and \ $H$ \ is an
ideal on \ $P_\kappa(\lambda)$ \ such that} \break $cof(H) < \hbox{{\bf
cov}}(\pmb{M}_{\kappa,\kappa}),$ \ {\sl and let \ $A\in H^+.$ \  Then 
there is \ $C\in
H^+\cap P(A)$ \ such that \ $[C]_\kappa^2 = [C]_\prec^2.$}
\bigskip
{\bf Proof.}\ For \ $\alpha<\kappa,$ \ set \ $A_\alpha = \{a\in A : 
\cup(a\cap\kappa)
= \alpha\}.$ \ By induction on \ $\alpha<\kappa,$ \ we define \ 
$c_k\in \{\phi\}\cup
A_\alpha$ \ for \ $k\in {}^\alpha 2$ \ as follows. Given \ $k\in 
{}^\alpha 2,$ \ set
$$e_k = \bigcup\{c_{k\upharpoonright\beta} : \beta\in k^{-1}(\{1\})\}$$
and
$$Z_k = \{a\in A_\alpha : e_k\subseteq a\}.$$
If \ $Z_k\not=\phi,$ \ let \ $c_k$ \ be an arbitrary member of \ 
$Z_k.$ \ Otherwise  let \
$c_k = \phi.$
\bigskip
Set \ $\nu = cof(H)$ \ and pick \ $B_\xi\in H$ \ for \ $\xi<\nu$ \ so 
that \ $H =
\displaystyle\bigcup_{\xi<\nu} P(B_\xi).$ \ Let \ $\xi<\nu.$ \ For \ 
$\alpha<\kappa,$
\ let \ $D_\xi^\alpha$\ be the set of all \  $s\in{}^{(\alpha+1)} 2$ 
\ such that (i)
\ $s(\alpha) = 1,$ \ and (ii) \ there is \  $a\in A_\alpha-B_\xi$ \ 
with the property that
$$(\forall\beta\in\alpha\cap s^{-1}(\{1\}))(\forall k\in{}^\beta 2)\quad
c_k\subseteq a.$$
Then let \ $D_\xi = \displaystyle\bigcup_{\alpha<\kappa} 
D_\xi^\alpha$ \ and \ $U_\xi
= \displaystyle\bigcup_{s\in D_\xi} O^\kappa_s.$ \ Let us prove that 
the open set \ $U_\xi$
\  is dense. Thus let \ $\gamma<\kappa$ \ and \ $p\in{}^\gamma 2.$ \ Pick
\ $a\in\displaystyle(\bigcup_{\gamma\leq\delta<\kappa} 
A_\delta)-B_\xi$ \ so that
$$(\forall\beta\in p^{-1}(\{1\}))(\forall k\in{}^\beta 2)\quad c_k\subseteq
a.$$
Set \ $\alpha = \cup(a\cap\kappa)$ \ and define \ $s\in{}^{(\alpha+1)}2$ \ by~:
\ $s\upharpoonright\gamma = p,\ s(\delta) = 0$ \ if \ 
$\gamma\leq\delta<\alpha,$ \ and
\ $s(\alpha) = 1.$ \ It is immediate that \ $s\in D_\xi^\alpha.$
\bigskip
Select \ $f\in\displaystyle\bigcap_{\xi<\nu} U_\xi.$ \ For each \ 
$\xi<\nu,$ \ there is
\ $s_\xi\in D_\xi$ \ such that \ $s_\xi\subset f.$ \ Let \ 
$\alpha_\xi<\kappa$ \ be such
that \ $s_\xi\in D_\xi^{\alpha_\xi}.$ \ Set \ $T = \{\alpha_\xi : 
\xi<\nu\}$ \ and
define \ $g\in{}^\kappa 2$ \ so that \ $g^{-1}(\{1\}) = T.$ \ For \ 
$\xi<\nu,$ \ set
$$d_\xi = \bigcup\{c_{g\upharpoonright\beta} : \beta\in T\cap\alpha_\xi\}$$
and
$$C_\xi = \{b\in A_{\alpha_\xi} : d_\xi\subseteq b\}.$$
Finally, let \ $C = \displaystyle\bigcup_{\xi<\nu}C_\xi.$
\bigskip
Let us verify that \ $C$ \ is as desired. It is clear that \ 
$C\subseteq A.$ \ Let
\ $\xi<\nu.$ \ There is \ $a_\xi\in A_{\alpha_\xi}-B_\xi$ \ such that
$$(\forall\beta\in\alpha_\xi\cap s_\xi^{-1}(\{1\}))(\forall k\in{}^\beta
2)\quad
c_k\subseteq a_\xi.$$
Put \ $k_\xi = g\upharpoonright\alpha_\xi.$ \ Then \ $a_\xi\in 
Z_{k_\xi}$ \ since
\ $s_\xi(\beta) = f(\beta) = 1$ \ for every \ $\beta\in 
T\cap\alpha_\xi.$ \ It follows
that \ $c_{k_\xi}\in Z_{k_\xi}.$ \ It is immediate that \ $Z_{k_\xi} = C_\xi.$
Thus we
have  shown that (a) $C_\xi-B_\xi\not=\phi$ \ for every \ $\xi<\nu,$ \ and (b)
$c_{g\upharpoonright\alpha_\xi}\in C_\xi$ \ for every \ $\xi<\nu.$ \ 
It follows from
(a) that \ $C\in H^+,$ \ and from (b) that \ $[C]_\kappa^2 = 
[C]_\prec^2$ \ since given
\ $\xi, \zeta<\nu$ \ with \ $\alpha_\xi<\alpha_\zeta,$ \ we have
$c_{g\upharpoonright\alpha_\xi}\subseteq b$ for every $b\in
C_\zeta.$\hfill$\carreblanc$
\bigskip\bigskip
{\bf QUESTION.} \ Is the assumption that \ $\kappa$ \ is inaccessible 
necessary in the
statement of Theorem 9.4 ?
\bigskip
{\bf Remark.} \ Suppose \ $\kappa$ \ is inaccessible. Then by 
Theorems 9.1, 9.2, 9.4, 5.4 and
4.7, \ $I_{\kappa,\lambda}$ \ is a weak \ $\chi$-point if and only if \
$\lambda^{<\kappa}<\hbox{\bf cov}(\pmb{M}_{\kappa,\kappa})$ \ if and only if \
$\{C:[C]_\kappa^2=[C]_\prec^2\}$ \ is dense in \ 
$(I_{\kappa,\lambda}^+,\subseteq)$.
\eject
\vskip 1,5cm
{\Bbf{10. $H^+ \hfltoutpetit{\kappa}{} (H^+,\alpha)^2$}}
\bigskip

{\bf Definition.} \ Let \ $H$ \ be an ideal on \ $P_\kappa(\lambda)$ 
\ and \ $\alpha$ \ an
ordinal. \
$H^+\hfltoutpetit{\kappa}{}(H^+,\alpha)^2$ \  means that given \  $F :
[P_\kappa(\lambda)]_\kappa^2 \hfltoutpetit{}{} 2$ \  and \ $A\in
H^+,$ \ there is  \
$B\subseteq A$  \ such that either \ $B\in H^+$ \ and \ $F$ \ is 
identically 0 on  \
$[B]_\kappa^2$ \  or \ $(B,\prec)$ \ has
order type \ $\alpha$ \ and \ $F$ \ is identically 1 on \ $[B]_\kappa^2.$
\bigskip
In this section we show that \ $H^+
\hfltoutpetit{\kappa}{}(H^+,\omega+1)^2$ \ for
every ideal \ $H$ \ on \ $P_\kappa(\lambda)$ \ with \ $cof(H)
<\hbox{{\bf cov}}(\pmb{M}_{\kappa,\kappa}).$
\bigskip\bigskip
{\bf Definition.}\ Suppose that \ $H$ \ is an ideal on \ 
$P_\kappa(\lambda), A\in H^+$ \ and
\ $F : \kappa\times P_\kappa(\lambda)\hfltoutpetit{}{} 2.$  Then $(H,A,F)$ \ is
0-good if there is \
$D\in H^+\cap P(A)$ \ such that \ $\{b\in D : F(\cup(a\cap\kappa),b) = 1\}\in
H$ \ for any  \ $a\in D.$
\bigskip
The following is straightforward.
\bigskip\bigskip
{\bf LEMMA 10.1.}\ {\sl Suppose that \ $(H,A,F)$ \ is 0-good, where \ 
$H$ \ is an ideal on
\ $P_\kappa(\lambda)$ \ which is both a weak \ $\pi$-point and a weak 
\ $\chi$-point, \
$A\in H^+$ \ and \ $F : \kappa\times 
P_\kappa(\lambda)\hfltoutpetit{}{}2.$ \ Then \
$F$ \ is identically 0 on \ $[C]_\kappa^2$ \ for some \ $C\in H^+\cap P(A).$}
\bigskip\bigskip
{\bf Definition.}\ Given an ideal \  $H$ \  on \ $P_\kappa(\lambda)$ 
\ and \ $B\in
H^+,$ \ let \ $M_{H,B}^{\hbox{d}}$ \ be the set of all \  $Q\subseteq 
H^+\cap P(B)$ \ such
that (i) any two distinct members of \ $Q$ \ are disjoint, and (ii) 
for every \ $A\in
H^+\cap P(B),$  \ there is \ $C\in Q$ \ with \ $A\cap C\in H^+.$
\bigskip\bigskip
{\bf LEMMA 10.2.}\ {\sl Suppose that \ $(H,A,F)$ \ is not 0-good, 
where \ $H$ \ is an ideal
on\
$P_\kappa(\lambda),\ A\in H^+$ \ and \ $F : \kappa\times
P_\kappa(\lambda)\hfltoutpetit{}{}2,$ \ and let \ $B\in H^+\cap 
P(A).$ \ Then there
exist}  \
$Q_B\in M_{H,B}^{\hbox{d}}$ \ {\sl and \
$\varphi_B : Q_B \hfltoutpetit{}{} B$ \ such that (i) 
$\varphi_B(D)\prec b$ \ and
\ $F(\cup(\varphi_B(D)\cap\kappa),b) = 1$ \ whenever \ $b\in D\in 
Q_B,$ \ and (ii)
$\cup(\varphi_B(D)\cap\kappa)\not=\cup(\varphi_B(D')\cap\kappa)$ \ for any two
distinct members \ $D$ \ and \ $D'$ \ of \ $Q_B.$}
\bigskip
{\bf Proof.} \ Set \ $T = \{\cup(a\cap\kappa) : a\in B\}$ \ and 
define \ $\psi :
T\times(H^+\cap P(B))\fle P(B)$ \ by \ $\psi(\alpha,C) = \{b\in C :
F(\alpha,b) = 1\}.$ \ Now using induction, define \ $\eta\leq\kappa$ 
\ and \ $\alpha_\delta\in T$ \
and \ $B_\delta\in H^+\cap P(B)$ \ for \ $\delta<\eta$ \ so that~:
\medskip
(0)\qquad If \ $\delta<\eta, \ B-(\displaystyle\bigcup_{\xi<\delta}B_\xi)\in
H^+,$
\medskip
$$\alpha_\delta = \hbox{ least } \alpha\in T \hbox{ such that }
\psi(\alpha,B-(\bigcup_{\xi<\delta}B_\xi))\in H^+$$
and \ $B_\delta =
\psi(\alpha_\delta,B-\displaystyle(\bigcup_{\xi<\delta}B_\xi)).$
\medskip
(1)\qquad If \ $\eta<\kappa, \
B-(\displaystyle\bigcup_{\xi<\eta}B_\xi)\in H.$
\medskip
Notice that if \  $\gamma<\delta<\eta,$ \ then
$$\psi(\alpha_\delta,B-(\bigcup_{\xi<\delta}B_\xi))\subseteq\psi(\alpha_\delta,
B-(\bigcup_{\zeta<\delta}B_\zeta))$$
and consequently \ $\alpha_\gamma\leq\alpha_\delta.$ \ In fact \
$\alpha_\gamma<\alpha_\delta$ \ as \
$\psi(\alpha_\gamma,B-(\displaystyle\bigcup_{\xi<\delta}B_\xi)) = 
\phi$ \ (since
\ $(B-\displaystyle\bigcup_{\xi<\delta}B_\xi))\cap B_\gamma = \phi$ \ 
and \ $B_\gamma
= \{b\in B-(\displaystyle\bigcup_{\zeta<\gamma}B_\zeta) : F(\alpha_\gamma,b)
= 1\}).$
\bigskip
We claim that \ $\{B_\delta : \delta<\eta\}\in M^
{\hbox{d}}_{H,B}.$ \ Suppose
otherwise. Then there exists \ $E\in H^+\cap P(B)$ \ such that \ 
$E\cap B_\xi\in
H$ \ for
every \ $\xi<\eta.$ \ Since
$$E-(\bigcup_{\xi<\delta}B_\xi) \in H^+\cap P(B-(\bigcup_{\xi<\delta}B_\xi))$$
for every \ $\delta<\kappa,$ \ we must have \ $\eta=\kappa.$ \ Set
$$\beta = \hbox{ least }\alpha\in T \hbox{ such that } \psi(\alpha,E)\in H^+.$$
Then for each \ $\delta<\kappa,$
$$\psi(\beta,E)-(\bigcup_{\xi<\delta}B_\xi)\in H^+\cap
P(\psi(B,B-(\bigcup_{\xi<\delta}B_\xi)))$$
and therefore \ $\beta\geq\alpha_\delta,$ \ which is a contradiction.
\bigskip
For each \ $\delta<\eta,$ \ pick \ $s_\delta\in B$ \ so that \ 
$\cup(s_\delta\cap\kappa) =
\alpha_\delta,$ \ and put
$$S_\delta = \{b\in B_\delta : s_\delta\cup(\alpha_\delta+2)\subseteq b\}.$$
Finally, set \ $Q_B = \{S_\delta : \delta<\eta\}$ \ and define \ 
$\varphi_B : Q_B
\hfltoutpetit{}{}
B$ \ by \ $\varphi_B(S_\delta) = s_\delta.$ \hfill$\carreblanc$
\bigskip\bigskip
{\bf LEMMA 10.3.} \ {\sl Suppose that \ $H$ \ is an ideal on \ 
$P_\kappa(\lambda)$ \ and
\ $A\in H^+.$ \ Suppose further that \ $C\in H^+\cap P(A)$ \ and } \ 
$Q_\alpha\in
M_{H,A}^{\hbox{d}}$ \ {\sl for \ $\alpha<\beta,$ \ where \ $\beta$ \ 
is a limit ordinal with
\ $0<\beta<\kappa.$ \ Then
$$\{a\in C : (\forall h\in \prod_{\alpha<\beta}Q_\alpha)\quad
a\notin\bigcap_{\alpha<\beta}h(\alpha)\}\in H.$$}
\bigskip
{\bf Proof.}\ It suffices to observe that for each
\ $a\in\displaystyle\bigcap_{\alpha<\beta}(C\cap(\cup Q_\alpha)),$ \ there is
\ $h\in\displaystyle\prod_{\alpha<\beta}Q_\alpha$ \ such that \break
$\displaystyle a\in\bigcap_{\alpha<\beta}h(\alpha).$~ \hfill$\carreblanc$
\bigskip\bigskip
{\bf LEMMA 10.4.}\ {\sl Suppose that \ $(H,A,F)$ \ is not 0-good, 
where \ $H$ \ is an ideal
on \ $P_\kappa(\lambda), \  A\in H^+$ \ and \ $F :
\kappa\times P_\kappa(\lambda)\hfltoutpetit{}{} 2.$ \ Then~:}
\medskip
{\parindent1,5cm
\item{\bf (i)} {\sl There is \ $C\subseteq A$ \ such that \ 
$(C,\prec)$ \ has order type \
$\omega+1$ \ and \ $F$ \ is identically 1 on \ $[C]_\kappa^2.$}
\medskip
\item{\bf (ii)} {\sl Suppose that \ $\hbox{\nomgot a}_H >\kappa$ \ 
and \ $\theta$ \ is
uncountable cardinal \ $<\kappa$ \ such that \ 
$\kappa\hfltoutpetit{}{} (\kappa,\theta)^2.$
\ Then there is \ $C\subseteq A$ \ such that \ $(C,\prec)$ \ has 
order type \ $\theta+1$ \
and \ $F$ \ is identically 1 on \ $[C]_\kappa^2.$}
\par }
\bigskip
{\bf Proof.} \ We prove (ii) and leave the proof of (i) to the reader.
By Corollary 19.7 in [EHM\'aR], we have that \ $\mu^\tau<\kappa$ \
whenever \ $\mu$ \ and \  $\tau$ \  are cardinals such that \ 
$\theta\leq\mu<\kappa$ \ and
\ $0<\tau<\theta.$ \ Using this and Lemmas 10.2 and 10.3, define
\ $R_\beta,Q_\beta\in\{W\in M_{H,A}^{\hbox{d}} : \vmid W<\kappa\}$ \ and
\ $\varphi_\beta : Q_\beta \hfltoutpetit{}{} A$ \ for \ $\beta<\theta$ by~:
\bigskip
(0)\qquad $R_0 = \{A\}$\ ;
\bigskip
(1)\qquad $Q_\beta = \displaystyle\bigcup_{B\in R_\beta}Q_B$\ ;
\bigskip
(2)\qquad $R_{\beta+1} = Q_\beta$ ;
\bigskip
(3)\qquad $R_\beta = H^+\cap\{\displaystyle\bigcap_{\alpha<\beta} h(\alpha) :
h\in\prod_{\alpha<\beta} Q_\alpha\}$ \ if \ $\beta$ \ is a limit 
ordinal \ $>0$ ;
\bigskip
(4)\qquad $\varphi_\beta = \displaystyle\bigcup_{B\in R_\beta}\varphi_B.$
\bigskip
Select \ $b\in \displaystyle\bigcap_{\beta<\theta}(\cup Q_\beta).$ \ There
must be \ $k\in\displaystyle\prod_{\beta<\theta} Q_\beta$ \ such that
\ $b\in\displaystyle\bigcap_{\beta<\theta} k(\beta).$ \ Then
$$C = \{\varphi_\beta(k(\beta)) : \beta<\theta\}\cup\{b\}$$
is as desired. \hfill$\carreblanc$
\bigskip\bigskip
{\bf THEOREM 10.5.} \ {\sl Suppose \ $\theta$ \ is
an infinite cardinal \ $<\kappa$ \ such that
\ $\kappa\hfltoutpetit{}{}(\kappa,\theta)^2.$ \ Then\break
\ $H^+\hfltoutpetit{\kappa}{} (H^+,\theta+1)^2$ \ for every ideal \ $H$ \ on
\ $P_\kappa(\lambda)$ \ such that } \ $cof(H)<\hbox{{\bf
cov}}(\pmb{M}_{\kappa,\kappa}).$
\bigskip
{\bf Proof.} \ Let \ $H$ \ be an ideal on \ $P_\kappa(\lambda)$ \ such that
\ $cof(H)<\hbox{{\bf cov}}(\pmb{M}_{\kappa,\kappa}).$ \ Then \ $H$ \  is a weak
\  $\chi$-point by
Theorem 9.1. Moreover, \ $H$ \ is a weak \ $\pi$-point and \ $\hbox{\nomgot
a}_H>\kappa$ \ by Theorem 8.4 since \ $\hbox{{\bf
cov}}(\pmb{M}_{\kappa,\kappa})\leq\hbox{\nomgot{d}}_{\kappa,\lambda}^\kappa$ 
\ by
Proposition 4.1. Hence, \ $H^+\hfltoutpetit{\kappa}{} (H^+,\theta+1)^2$ \
by Lemmas 10.1 and 10.4.\hfill$\carreblanc$

\vskip 1,5cm
{\Bbf{11. \pmb{$H^+ \hfltoutpetit{\kappa}{\kappa} \ (H^+,\alpha)^2$}}}
\bigskip
{\bf Definition.} \ For \ $A\subseteq P_\kappa(\lambda),$ \ let
$$[A]_{\kappa,\kappa}^2 = \{(\cup(a\cap\kappa),\cup(b\cap\kappa)) : 
a,b\in A \hbox{ and }
\cup(a\cap\kappa) < \cup(b\cap\kappa)\}$$
{\bf Remark.} \
$[P_\kappa(\lambda)]_{\kappa,\kappa}^2 = [\kappa]^2.$
\bigskip\bigskip
{\bf Definition.} \ Let \ $H$ \ be an ideal on \ $P_\kappa(\lambda)$ 
\ and \ $\alpha$ \ an
ordinal. \ $H^+\hfltoutpetit{\kappa}{\kappa} \ (H^+,\alpha)^2$ \ 
means that given \ $F :
[P_\kappa(\lambda)]^2_{\kappa,\kappa}\hfltoutpetit{}{} 2$ \ and \ 
$A\in H^+,$ \ there is \
$B\subseteq A$ \ such that either \ $B\in H^+$ \ and \ $F$ \ is 
identically 0 on \
$[B]_{\kappa,\kappa}^2,$ \ or \ $(B,\prec)$ \ has order type \ 
$\alpha$ \ and \ $F$ \ is
identically 1 on \ $[B]_{\kappa,\kappa}^2.$
\bigskip
We will show that \ $H^+ \hfltoutpetit{\kappa}{\kappa} \ 
(H^+,\omega+1)^2$ \ for every
ideal \
$H$ \ on \ $P_\kappa(\lambda)$ \ such that \ $cof(H) < \hbox{\bf 
non}_\kappa$(weakly
selective).
\bigskip\bigskip
{\bf Definition.} \ For an ideal \ $H$ \ on \ $P_\kappa(\lambda), \ 
J_H = \ \{B\subseteq
\kappa : U_B\in H\},$ \ where
$$U_B = \{a\in P_\kappa(\lambda) : \cup(a\cap\kappa)\in B\}.$$
{\bf LEMMA 11.1.} \ {\sl Let \ $H$ \ be an ideal on  \ 
$P_\kappa(\lambda).$  \ Then \ $J_H$
\ is an ideal on \ $\kappa.$ \ Moreover, \ $cof(J_H) \leq cof(H).$}
\bigskip
{\bf Proof.} \ It is simple to see that (a) $U_\kappa = P_\kappa(\lambda),$ (b)
$U_{\cup\hbox{\nomgot 
B}}\subseteq\displaystyle\bigcup_{B\in{\hbox{\nomgot B}}} U_B$ \ for \
$\hbox{\nomgot B}\subseteq P(\kappa),$ (c) $U_C\subseteq U_B$ \ if \ 
$C\subseteq
B\subseteq  K,$ \ and (d)
$U_B\in I_{\kappa,\lambda}$ \ for every \ $B\subseteq\kappa$ \ with \ 
$\vmid B=1.$ \ The
first assertion immediately follows.
\bigskip
For \ $C\subseteq P_\kappa(\lambda),$ \ let \ $Y_C$ \ be the set of 
all \ $\delta\in\kappa$
\ such that
$$\{a\in P_\kappa(\lambda) : \cup(a\cap\kappa) = \delta\}\subseteq C.$$
If \ $C\in H,$ \ then \ $Y_C\in J_H$ \ since \ $U_{Y_C}\subseteq C.$ 
\ Moreover if \
$B\subseteq\kappa$ \ and \ $U_B\subseteq C\subseteq 
P_\kappa(\lambda),$ \ then \
$B\subseteq Y_C.$ \ Hence \ $cof(J_H)\leq cof(H).$ \hfill\carreblanc
\bigskip
{\bf Remark.} \ Let \ $H$ \ be an ideal on \ $P_\kappa(\lambda).$ \ Then
$$\{\cup(a\cap\kappa) : a\in A\}\in (J_{H\mid A})^+$$
for every \ $A\in H^+.$
\bigskip
The following is readily checked.
\bigskip\bigskip
{\bf LEMMA 11.2.} \ {\sl Given an ideal \ $H$ \ on \ 
$P_\kappa(\lambda),$ \ the following
are equivalent~:}
\medskip
{\parindent=1,5cm
\item{\bf (i)} {\sl $J_H$ \ is a local \ $Q$-point.}
\medskip
\item{\bf (ii)} {\sl For every \ $g\in{}^\kappa\kappa,$ \ there is \ 
$B\in H^+$ \ such that
\ $g(\cup(a\cap\kappa))<\cup(b\cap\kappa)$ \ for all \ $a,b\in B$ \ with \
$\cup(a\cap\kappa)<\cup(b\cap\kappa).$}
\par}
\bigskip
Suppose \ $\kappa$ \ is a limit cardinal. If \ $\kappa^+<\hbox{\bf 
non}_\kappa$(weak
$Q$-point), then by Lemma 11.1 \ $J_{I_{\kappa,\kappa^+}\mid A}$ \ is 
a local \ $Q$-point
for every \ $A\in I_{\kappa,\kappa^+}^+.$ \ The following shows that 
this implication can
be reversed.
\bigskip\bigskip
{\bf PROPOSITION 11.3.} \ {\sl Suppose that \ $\kappa$ \ is a limit 
cardinal and \
$J_{I_{\kappa,\lambda}\mid A}$ \ is a local \ $Q$-point for every \ $A\in
I_{\kappa,\lambda}^+.$ \ Then }\ $\sigma<\overline{\hbox{\bf 
non}}_\kappa$(weak $Q$-point)
{\sl for every
\ $\sigma\in{\cal K}(\kappa,\lambda).$
\par}
\bigskip
{\bf Proof.} \ Suppose that \ $J$ \ is an ideal on \ $\kappa$ \ and \ 
$T\subseteq
P_\kappa(\lambda-\kappa)$ \ is such that \ 
$\overline{cof}(J)\leq\vmid T$ \ and \
$\vmid{T\cap P(a)}<\kappa$ \ for every \ $a\in P_\kappa(\lambda).$ \ 
Select \ $B_d\in J$ \
for \ $d\in T$ \ so that for every \ $D\in J,$ \ there is \ $u\in 
P_\kappa(T)-\{\phi\}$ \
with \ $D\subseteq\displaystyle\bigcup_{d\in u}B_d.$ \ Let \ $A$ \ be 
the set of all \ $a\in
P_\kappa(\lambda)$ \ such that \ $\cup(a\cap\kappa)\notin B_d$ \ for 
every \ $d\in T\cap
P(a-\kappa).$ \ It is simple to see that \ $A\in 
I_{\kappa,\lambda}^+.$ \ Now fix \
$g\in{}^\kappa\kappa.$ \ By Lemma 11.2, there is \ $C\in 
(I_{\kappa,\lambda}\mid A)^+$ \
such that \ $g(\cup(a\cap\kappa))<\cup(b\cap\kappa)$ \ for all \ 
$a,b\in C$ \ with \
$\cup(a\cap\kappa)<\cup(b\cap\kappa).$ \ Set
$$D = \{\cup(a\cap\kappa) : a\in C\cap A\}.$$
Then \ $D\in J^+.$ \ Moreover \ $g(\alpha)<\beta$ \ for all \ 
$\alpha,\beta\in D$ \ with \
$\alpha<\beta.$ \ Hence \ $J$ \ is a local $Q$-point. \hfill\carreblanc
\bigskip\bigskip
{\bf THEOREM 11.4.} \ {\sl Suppose that \ $\theta$ \ is an infinite 
cardinal \ $<\kappa$
\ such that \ $\kappa\hfltoutpetit{}{} (\kappa,\theta)^2,$ \ and \ 
$H$ \ is an ideal on \
$P_\kappa(\lambda)$ \ with }\ $cof(H)<\hbox{\bf non}_\kappa$(weakly 
selective). {\sl Then \
$H^+
\hfltoutpetit{\kappa}{\kappa} (H^+,\theta+1)^2.$
  }
\bigskip
{\bf Proof.} \ Fix \ $G :\kappa\times\kappa\hfltoutpetit{}{} 2$ \ and 
\ $A\in H^+$. \
Define \
$F : \kappa\times P_\kappa(\lambda)\hfltoutpetit{}{} 2$ by \ $F(\alpha,b) =
G(\alpha,\cup(b\cap\kappa)).$
\bigskip
First suppose \ $(H,A,F)$ \ is 0-good. Pick \ $D\in
H^+\cap P(A)$ \ so that
$$\{b\in D : F(\cup(a\cap\kappa),b) = 1\}\in H$$
for any \ $a\in D.$ \ Set \ $B_\alpha = \{\delta<\kappa : 
G(\alpha,\delta) = 1\}$ \ for \
$\alpha<\kappa.$ \ Then \ $B_{\cup(a\cap\kappa)}\in J_{H\mid D}$ \ 
for every \ $a\in D$ \
since
$$D\cap U_{B_{\cup(a\cap\kappa)}} = \{b\in D : 
G(\cup(a\cap\kappa),\cup(b\cap\kappa)) =
1\}= \{b\in D : F(\cup(a\cap\kappa),b)=1\}.$$
By Lemma 11.1 \ $cof(J_{H\mid D}) < \hbox{\bf non}_\kappa$(weak 
$P$-point) so there is \
$G\in (J_{H\mid D})^+$ \ such that \ $\vmid{G\cap 
B_{\cup(a\cap\kappa)}}<\kappa$ \ for
every \ $a\in D.$ \ Notice that \ $D\cap U_G\in H^+.$ \ Select \ 
$g\in{}^\kappa\kappa$ \ so
that \ $\cup(b\cap\kappa)\notin B_{\cup(a\cap\kappa)}$ \ for all \ 
$a,b\in D\cap U_G$ \
such that \ $g(\cup(a\cap\kappa))<\cup(b\cap\kappa).$ \ By Lemma 11.1
$$cof(J_{H\mid(D\cap U_G)}) <\hbox{\bf non}_\kappa\hbox{(weak } 
Q\hbox{-point)}$$
and hence by Lemma 11.2 there is \ $R\in(H\mid(D\cap U_G))^+$ \ such that \
$g(\cup(a\cap\kappa))<\cup(b\cap\kappa)$ \ for all \ $a,b\in R$ \ with \
$\cup(a\cap\kappa)<\cup(b\cap\kappa).$ \ Then \ $R\cap D\cap U_G\in 
H^+\cap P(A)$ \ and
moreover \ $F$ \ is identically 0 on \ $[R\cap D\cap U_G]_{\kappa,\kappa}^2.$
\bigskip
Finally, suppose \ $(H,A,F)$ \ is not 0-good. Since \ $\hbox{\nomgot 
a}_H>\kappa$ \ by
Theorems 2.2 and 8.4, there is by Lemma 10.4 \ $C\subseteq A$ \ such 
that \ $(C,\prec)$ \
has order type \ $\theta+1$ \ and \ $F$ \ is identically 1 on \ 
$[C]_\kappa^2.$ \ It is
immediate that \ $G$ \ is constantly 1 on \ $[C]_{\kappa,\kappa}^2.$ 
\hfill\carreblanc
\bigskip
{\bf Remark.} \ Suppose \ $\kappa$ \ is a successor cardinal. Then by 
Theorem 11.4 \
$\kappa^+<\hbox{\nomgot d}_\kappa$ \ implies that \break $I_{\kappa,\kappa^+}^+
\hfltoutpetit{\kappa}{\kappa} (I_{\kappa,\kappa^+}^+,\theta+1)^2$ \ for 
every cardinal \
$\theta\geq 2$ \ such that \ $\kappa\hfltoutpetit{}{} 
(\kappa,\theta)^2.$ \ Conversely, it
will be shown in the next section that \ 
$I_{\kappa,\kappa^+}^+\hfltoutpetit{\kappa}{\kappa}
(I_{\kappa,\kappa^+}^+,3)^2$ \ implies that \ $\kappa^+<\hbox{\nomgot 
d}_\kappa.$
\vskip 1,5cm
{\Bbf{12. \pmb{$H^+ \hfltoutpetit{\kappa}{\kappa} \ (H^+;\alpha)^2$}}}
\bigskip
{\bf Definition.} \ Given an ideal \ $H$ \ on \ $P_\kappa(\lambda)$ \ 
and an ordinal \
$\alpha, H^+\hfltoutpetit{\kappa}{\kappa} (H^+;\alpha)^2$ \ means 
that for all \break $F :
[P_\kappa(\lambda)]_{\kappa,\kappa}^2\hfltoutpetit{}{} 2$ \ and \ 
$A\in H^+,$ \ there is \
$B\subseteq A$ \ such that either \ $B\in H^+$ \ and \ $F$ \ is 
identically 0 on \
$[B]_{\kappa,\kappa}^2,$ \ or \ $\{\cup(a\cap\kappa) : a\in B\}$ \ 
has order type \
$\alpha$ \ and \ $F$ \ is identically 1 on \ $[B]_{\kappa,\kappa}^2.$
\bigskip
{\bf Remark.} \ $H^+\hfltoutpetit{\kappa}{} (H^+,\alpha)^2\Rightarrow
H^+\hfltoutpetit{\kappa}{\kappa}(H^+,\alpha)^2\Rightarrow 
H^+\hfltoutpetit{\kappa}{\kappa}
(H^+;\alpha)^2\Rightarrow \kappa\hfltoutpetit{}{} (\kappa,\alpha)^2.$ \
\bigskip
  We will prove that \
$I_{\kappa,\kappa}^+\hfltoutpetit{\kappa}{\kappa} 
(I_{\kappa,\kappa^+}^+;\alpha)^2$ \ if and
only if \ $\kappa^+<\hbox{\bf non}_\kappa(J^+\hfltoutpetit{}{} 
(J^+,\alpha)^2).$
\bigskip\bigskip
{\bf THEOREM 12.1.} \ {\sl Suppose that \ $3\leq\alpha\leq\kappa$ \ 
and \ $H$ \ is an
ideal on \ $P_\kappa(\lambda)$ \ such that \break $cof(H)<\hbox{\bf
non}_\kappa(J^+\hfltoutpetit{}{}(J^+,\alpha)^2).$ \ Then \ 
$H^+\hfltoutpetit{\kappa}{\kappa}
(H^+;\alpha)^2.$}
\bigskip
{\bf Proof.} \ By Lemma 11.1, \ $(J_{H\mid A})^+ \hfltoutpetit{}{}((J_{H\mid
A})^+,\alpha)^2$ \ for every \ $A\in H^+.$ \ The desired conclusion 
easily follows.~
\hfill\carreblanc
\bigskip\bigskip
{\bf THEOREM 12.2.} \ {\sl Suppose that \ $3\leq\alpha\leq\kappa$ \ and $ \
I_{\kappa,\lambda}^+\hfltoutpetit{\kappa}{\kappa} 
(I_{\kappa,\lambda}^+;\alpha)^2.$ \ Then \
$\sigma<\overline{\hbox{\bf non}}_\kappa (J^+\hfltoutpetit{}{} 
(J^+,\alpha)^2)$ \ for every
\
$\sigma\in{\cal K}(\kappa,\lambda).$
}
\bigskip
{\bf Proof.} \  The proof is an easy modification of that of Proposition
11.3. \hfill\carreblanc
\bigskip
{\bf Remark.} \ Suppose that \ $\kappa$ \ is inaccessible and \ 
$3\leq\alpha\leq\kappa.$ \
Then by Theorems 12.1 and 12.2, \ $I_{\kappa,\lambda}^+ 
\hfltoutpetit{\kappa}{\kappa}
(I_{\kappa,\lambda}^+;\alpha)^2$ \ if and only if \ 
$\lambda^{<\kappa}<\hbox{\bf
non}_\kappa(J^+\hfltoutpetit{}{} (J^+,\alpha)^2).$
\bigskip
Let us finally observe that for \ $3\leq\alpha\leq\kappa,$ \ there 
always exists an ideal \
$H$ \ on  \ $P_\kappa(\lambda)$ \ of the least possible cofinality such that \
$H^+\hfltoutpetit{\kappa}{\kappa}\kern-3mm/\kern2mm (H^+;\alpha)^2$~:
\bigskip\bigskip
{\bf PROPOSITION 12.3.} \ {\sl Given \ $3\leq\alpha\leq\kappa,$ \ 
there is an ideal \ $H$ \
on \ $P_\kappa(\lambda)$ \ such that} (a) {\sl 
$H^+\hfltoutpetit{\kappa}{\kappa}\kern-3mm/\kern2mm
(H^+;\alpha)^2,$} \break
(b) \  {\sl $cof(H) = u(\kappa,\lambda)\cdot\hbox{\bf 
non}_\kappa(J^+\hfltoutpetit{}{}
(J^+,\alpha)^2),$ \ and }
(c) {\sl $\overline{cof}(H) \leq \lambda\cdot\overline{\bf
non}_\kappa(J^+\hfltoutpetit{}{} (J^+,\alpha)^2).$
  }
\bigskip
{\bf Proof.} \ Argue as for Lemma 5.1 of [M2]. \hfill\carreblanc

\eject
{\Bbf{13. \pmb{$H^+ \hfltoutpetit{\kappa}{} \ (H^+)^2$}}}
\bigskip
{\bf Definition.} \ Given an ideal \ $H$ \ on \ $P_\kappa(\lambda), \
H^+\hfltoutpetit{\kappa}{} (H^+)^2$ \ (respectively, \ 
$H^+\hfltoutpetit{\kappa}{\kappa}
(H^+)^2)$ \ means that for all \ $F : 
[P_\kappa(\lambda)]_\kappa^2\hfltoutpetit{}{} 2$ \
(respectively, \ $F : 
[P_\kappa(\lambda)]_{\kappa,\kappa}^2\hfltoutpetit{}{} 2)$ \ and \
$A\in H^+,$ \ there is \
$B\in H^+\cap P(A)$ \ such that \ $F$ \ is constant on \ 
$[B]_\kappa^2$ \ (respectively, \
$[B]_{\kappa,\kappa}^2).$
\bigskip\bigskip
{\bf THEOREM 13.1.} \ {\sl Suppose \ $\kappa$ \ is
weakly compact. Then $H^+ \hfltoutpetit{\kappa}{} (H^+)^2$ for every 
ideal $H$ on
$P_\kappa(\lambda)$ such that} $cof(H)<\hbox{{\bf
cov}}(\pmb{M}_{\kappa,\kappa}).$
\bigskip
{\bf Proof.}\ Suppose that \ $H$ \ is an ideal on \ $P_\kappa(\lambda)$ \ with
\ $cof(H)< \hbox{{\bf cov}}(\pmb{M}_{\kappa,\kappa}), F : \kappa\times
P_\kappa(\lambda)\hfltoutpetit{}{} 2$ \ and \ $A\in H^+.$ \ Then \  $cof(H)
<\hbox{\nomgot{d}}_{\kappa,\lambda}^\kappa$ by Proposition 4.1 and
therefore by a result of [M5] there are \ $B\in H^+\cap P(A)$ \ and \ 
$i<2$ \ such that
$$\{b\in B : F(\cup(a\cap\kappa),b)\not= i\}\in I_{\kappa,\lambda}$$
  for every \
$a\in B.$ \ Since \ $H$ \ is a weak \ $\chi$-point by Theorem 9.1, 
there is \ $C\in
H^+\cap P(B)$ \ such that \ $F$ \ takes the constant value \ $i$ \ on
\ $[C]_\kappa^2.$\hfill$\carreblanc$
\bigskip
{\bf Remark.} \ It follows from Theorem 6.5 (ii) and Theorem 15.1 
(below) that if \ $\kappa$
\ is weakly compact, then \ $H^+\hfltoutpetit{\kappa}{\kappa} 
(H^+)^2$ \ for every ideal \
$H$
\ on \ $P_\kappa(\lambda)$ \ such that \ $cof(H) <\hbox{\bf 
non}_\kappa$(weakly selective).
\bigskip\bigskip
{\bf COROLLARY 13.2.} {\sl The following are equivalent~:}
\medskip
{\parindent=1,5cm
\item{\bf (i)} $\kappa$ \ is weakly compact and \ $\lambda^{<\kappa}<\hbox{\bf
cov}(\pmb{M}_{\kappa,\kappa}).$
\medskip
\item{\bf (ii)} $I_{\kappa,\lambda}^+ \hfltoutpetit{\kappa}{} 
(I_{\kappa,\lambda}^+)^2.$
\medskip
\item{\bf (iii)} $I_{\kappa,\lambda}^+\hfltoutpetit{\kappa}{\kappa}
(I_{\kappa,\lambda}^+;\kappa)^2.$
\par}
\bigskip
{\bf Proof.}
(i) \ \ $\fle$ (ii) : By Theorem 13.1.

(ii) \ $\fle$ (iii) : Trivial.

(iii) $\fle$ (i) : By Theorems 12.2, 6.5 (i), 6.1 (iii), 5.4 and 4.7. 
\hfill\carreblanc

\vskip1,5cm
{\Bbf{14. \pmb{$H^+ \hfltoutpetit{\kappa}{} \ [H^+]_\rho^2$}}}
\bigskip
{\bf Definition.} \ Given a cardinal \ $\rho$ \ with \ 
$2\leq\rho\leq\lambda^{<\kappa}$ \
and an ideal \ $H$ \ on \ $P_\kappa(\lambda), 
H^+\hfltoutpetit{\kappa}{} [H^+]_\rho^2$ \
means that for all \ $F : [P_\kappa(\lambda)]_\kappa^2 
\hfltoutpetit{}{}\rho$ \ and \ $A\in
H^+,$ \  there is \ $B\in H^+\cap P(A)$ such that \ $F''[B]_\kappa^2\not=\rho.$

\bigskip\bigskip
{\bf THEOREM 14.1.} \ {\sl Suppose that \ $\kappa$ \ is a limit 
cardinal and \ $H$ \ is
an ideal on \ $P_\kappa(\lambda)$ \ such that}
\ $cof(H)<\hbox{{\bf cov}}(\pmb{M}_{\kappa,\kappa}).$ \ {\sl Then
\ $H^+\hfltoutpetit{\kappa}{} [H^+]_{\kappa^+}^2.$}
\bigskip
{\bf Proof.} \ Fix \ $F : \kappa\times
P_\kappa(\lambda)\hfltoutpetit{}{}\kappa^+$ \ and \ $A\in H^+.$ \ 
Since \ $cof(H)
<\hbox{\nomgot d}_{\kappa,\lambda}^\kappa$ \ by Proposition 4.1,
there are \ $B\in H^+\cap P(A)$ \ and \ $\xi\in\kappa^+$ \ such that
\ $\{b\in B : F(\cup(a\cap\kappa),b) = \xi\}\in I_{\kappa,\lambda}$ \ for every
\ $a\in B$ \ ([M5]). Now \ $H$ \ is a weak \ $\chi$-point by Theorem 9.1 and so
\ $\xi\notin F''[C]_\kappa^2$ \ for some \ $C\in H^+\cap 
P(B).$\hfill$\carreblanc$
\bigskip
Let us now show that \ $I_{\kappa,\lambda}^+ 
\hfltoutpetit{\kappa}{}\kern-3mm/\kern2mm
[I_{\kappa,\lambda}^+]_\lambda^2$ \ if \ 
$\lambda\geq\overline{\hbox{\nomgot d}}_\kappa.$ \
We will need some definitions.
\bigskip\bigskip
{\bf Definition.} \ Given \ 
$f\in\displaystyle\prod_{\alpha\in\kappa}(\kappa-\alpha),$ \ we
define \ $\widetilde f\in{}^\kappa\kappa$ \ by stipulating that
\medskip
(i) $\widetilde f(0) = 0$ ;
\medskip
(ii) $\widetilde f(\xi+1) = f(\widetilde f(\xi))+1$ ;
\medskip
(iii) $\widetilde f(\xi) = \displaystyle\bigcup_{\zeta<\xi}\widetilde f(\zeta)$
\ if \ $\xi$ \ is a limit ordinal \ $>0.$

\bigskip
{\bf Remark.} \ $\widetilde f$ \ is a strictly increasing function.
\bigskip
{\bf Remark.} \ If \ $g\in{}^\kappa\kappa$ \ is a strictly increasing 
function such that \
$g(\alpha)\leq f(\alpha)$ \ for all \ $\alpha<\kappa,$ \ then \ $g(\widetilde
f(\xi))\in[\widetilde f(\xi),\widetilde f(\xi+1))$ \ for every \ $\xi<\kappa.$
\bigskip\bigskip
{\bf Definition.} \ Given \ 
$f\in\displaystyle\prod_{\alpha\in\kappa}(\kappa-\alpha)$ \ and
a cardinal \ $\tau\in(0,\kappa),$ \ we define \ $c_{f,\tau} : \widetilde
f(\tau)\hfltoutpetit{}{}\tau$ \ by stipulating that \ $c_{f,\tau}$ \ 
takes the constant
value \
$\xi$ \ on \ $[\widetilde f(\xi),\widetilde f(\xi+1)).$
\bigskip\bigskip
{\bf Definition.} \ Suppose that \ $T\subseteq 
P_\kappa(\lambda-\kappa)$ \ is such that (a)
$\vmid T\geq\overline{\hbox{\nomgot d}}_\kappa,$ \ and (b) 
$\vmid{T\cap P(a)}<\kappa$ \ for
every \ $a\in P_\kappa(\lambda).$

Let \ $\psi_T : T\hfltoutpetit{}{}{}^\kappa\kappa$ \ be such that given \
$g\in{}^\kappa\kappa,$ \ there is \ $u\in P_\kappa(T)-\{\phi\}$ \ such that
$$g(\alpha)\leq\bigcup_{d\in u}(\psi_T(d))(\alpha)$$
for all \ $\alpha<\kappa.$

For \ $e\in P_\kappa(\lambda-\kappa),$ \ let \ $\tau_{T,e} = 
\vmid{T\cap P(e)}$ \ and
select a bijection \ $k_{T,e} : \tau_{T,e}\hfltoutpetit{}{} T\cap P(e).$

Also, define \ $f_{T,e}\in{}^\kappa\kappa$ \ by
$$f_{T,e}(\alpha) = \hbox{max}(\alpha,\bigcup_{d\in T\cap 
P(e)}(\psi_T(d))(\alpha)).$$
Finally, let \ $A_T$ \ be the set of all \ $a\in P_\kappa(\lambda)$ \ 
such that (i) $T\cap
P(a-\kappa)\not=\phi,$ \ and (ii) $\cup(a\cap\kappa)\geq\widetilde
f_{T,a-\kappa}(\tau_{T,a-\kappa}).$
\bigskip
{\bf Remark.} \  $A_T\in I_{\kappa,\lambda}^+.$
\bigskip\bigskip
{\bf THEOREM 14.2.} \ {\sl Suppose that \ $\rho\in{\cal 
K}(\kappa,\lambda)$ \ and \
$\rho\geq\overline{\hbox{\nomgot d}}_\kappa.$ \ Then \
$I_{\kappa,\lambda}^+\hfltoutpetit{\kappa}{}\kern-3mm/\kern2mm
[I_{\kappa,\lambda}^+]_\rho^2.$  }
\bigskip
{\bf Proof.} \ Select \ $T\subseteq P_\kappa(\lambda-\kappa)$ \ so 
that \ $\vmid T = \rho$
\ and \ $\vmid{T\cap P(a)}<\kappa$ \ for every \ $a\in 
P_\kappa(\lambda).$ \ We define a
partial function \ $F$ \ from \ $\kappa\times A_T$ \ to \ $T$ \ by 
stipulating that
$$F(\beta,a) = k_{T,a-\kappa}(c_{f_{T,a-\kappa},\tau_{T,a-\kappa}}(\beta))$$
if \ $a\in A_T$ \ and \ $\beta<\widetilde f_{T,a-\kappa}(\tau_{T,a-\kappa}).$
\bigskip
Now fix \ $B\in I_{\kappa,\lambda}^+\cap P(A_T)$ \ and \ $x\in T.$ \ Let \
$g\in{}^\kappa\kappa$ \ be the increasing enumeration of the elements 
of the set \
$\{\cup(b\cap\kappa) : b\in B\}.$ \ Select \ $u\in 
P_\kappa(T)-\{\phi\}$ \ so that \
$g(\alpha)\leq\displaystyle\bigcup_{d\in u}(\psi_T(d))(\alpha)$ \ for 
all \ $\alpha<\kappa.$
\ Now pick
\ $a\in B$ \ so that \ $x\cup(\cup u)\subseteq a.$ \ Notice that \ 
$g(\alpha)\leq
f_{T,a-\kappa}(\alpha)$ \ for every \ $\alpha\in\kappa.$ \ Let \ 
$\xi\in\tau_{T,a-\kappa}$
\ be such that \ $k_{T,a-\kappa}(\xi) = x.$ \ Then
$$\widetilde f_{T,a-\kappa}(\xi)\leq g(\widetilde 
f_{T,a-\kappa}(\xi)) < \widetilde
f_{T,a-\kappa}(\xi+1)\leq\widetilde 
f_{T,a-\kappa}(\tau_{T,a-\kappa})\leq\cup(a\cap\kappa).$$
Moreover,
$$F(g(\widetilde f_{T,a-\kappa}(\xi)),a) = k_{T,a-\kappa}(\xi) = x.$$
since
$$c_{f_{T,a-\kappa},\tau_{T,a-\kappa}}(g(\widetilde 
f_{T,a-\kappa}(\xi))) = \xi$$
  \remonte\hfill\carreblanc

\vskip 1,5cm
{\Bbf{15. \pmb{$H^+ \hfltoutpetit{\kappa}{\kappa} \ [H^+]_\rho^2$}}}
\bigskip
{\bf Definition.} \ Given a cardinal \ $\rho\in[2,\kappa]$ \ and an 
ideal \ $H$ \ on \
$P_\kappa(\lambda),\ H^+\hfltoutpetit{\kappa}{\kappa} [H^+]_\rho^2$ \ 
means that for all \
$F : [P_\kappa(\lambda)]_{\kappa,\kappa}^2 \hfltoutpetit{}{}\rho$ \ 
and \ $A\in H^+,$ \
there is \
$B\in H^+\cap P(A)$ \ such that \ $F''[B]_{\kappa,\kappa}^2\not=\rho.$
\bigskip
{\bf Remark.} \ $\kappa\hfltoutpetit{}{}\kern-3mm/\kern2mm 
[\kappa]_\rho^2 \Rightarrow
H^+\hfltoutpetit{\kappa}{\kappa}\kern-3mm/\kern2mm [H^+]_\rho^2 \Rightarrow
H^+\hfltoutpetit{\kappa}{}\kern-3mm/\kern2mm [H^+]_\rho^2.$
\bigskip
The following result shows that \ 
$I_{\kappa,\kappa^+}^+\hfltoutpetit{\kappa}{\kappa}
[I_{\kappa,\kappa^+}^+]_\rho^2$ \ if and only if \ 
$\kappa^+<\hbox{\bf non}_\kappa
(J^+\hfltoutpetit{}{} [J^+]_\rho^2).$
\bigskip\bigskip
{\bf THEOREM 15.1.} \ {\sl Let \ $\rho$ \ be a cardinal with \ 
$2\leq\rho\leq\kappa.$ \
Then~:}
\medskip
{\parindent=1,5cm
\item{\bf (i)} $H^+\hfltoutpetit{\kappa}{\kappa} [H^+]_\rho^2$ \ {\sl 
for every ideal \ $H$
\ on
\ $P_\kappa(\lambda)$ \ such that \
$cof(H) <\hbox{\bf non}_\kappa(J^+\hfltoutpetit{}{} [J^+]_\rho^2).$}
\medskip
\item{\bf (ii)} {\sl If \ $I_{\kappa,\lambda}^+\hfltoutpetit{\kappa}{\kappa}
[I_{\kappa,\lambda}^+]_\rho^2,$ \ then \
$\sigma<\overline{\hbox{\bf non}}_\kappa(J^+\hfltoutpetit{}{} 
[J^+]_\rho^2)$ \ for every \
$\sigma\in{\cal K}(\kappa,\lambda).$}
\par}
\bigskip
{\bf Proof.}
(i) : Use Lemma 11.1.

(ii) : Argue as for Proposition 11.3. \hfill\carreblanc

\bigskip
{\bf Remark.} \ Thus assuming \ $\kappa$ \ is inaccessible, \ 
$I_{\kappa,\lambda}^+
\hfltoutpetit{\kappa}{\kappa} [I_{\kappa,\lambda}^+]_\rho^2$ \ if and only if \
$\lambda^{<\kappa}<\hbox{\bf non}_\kappa(J^+\hfltoutpetit{}{} [J^+]_\rho^2).$
\bigskip
Finally, we show that if \ $\lambda\geq\overline{\hbox{\nomgot 
d}}_\kappa$ \ and \ $\kappa$
\ is a limit cardinal such that \ $2^{<\kappa} = \kappa,$ \ then \
$I_{\kappa,\lambda}^+\hfltoutpetit{\kappa}{\kappa}\kern-3mm/\kern2mm
[I_{\kappa,\lambda}^+]_\kappa^2.$
\bigskip\bigskip
{\bf THEOREM 15.2.} \ {\sl Suppose that (a) $\kappa$ \ is a limit 
cardinal such that \
$2^{<\kappa}=\kappa,$ \ and (b) either \ 
$\lambda>\overline{\hbox{\nomgot d}}_\kappa,$ \ or
\
$\overline{\hbox{\nomgot d}}_\kappa\in{\cal K}(\kappa,\lambda).$ \ Then \
$I_{\kappa,\lambda}^+ \hfltoutpetit{\kappa}{\kappa}\kern-3mm/\kern2mm
[I_{\kappa,\lambda}^+]_\kappa^2.$
  }
\bigskip
{\bf Proof.} \ Select \ $T\subseteq P_\kappa(\lambda-\kappa)$ \ so 
that \ $\vmid T =
\lambda\cdot\overline{\hbox{\nomgot d}}_\kappa$ \ and \ $\vmid{T\cap 
P(a)}<\kappa$ \ for
every \ $a\in P_\kappa(\lambda).$ \ Also, select  \ $\chi :
\kappa\hfltoutpetit{}{}\displaystyle\bigcup_{\gamma<\kappa} 
{}^\gamma\kappa$ \ so that \
$\vmid{\chi^{-1}(\{z\})}=\kappa$ \ for every \ $z\in
\displaystyle\bigcup_{\gamma<\kappa}{}^\gamma\kappa.$
\ Now let
\
$A$
\ be the set of all \ $a\in A_T$ \ such that
$$\chi(\cup(a\cap\kappa)) = c_{f_{T,a-\kappa},\tau_{T,a-\kappa}}.$$
Notice that \ $A\in I_{\kappa,\lambda}^+.$ \ We define a partial 
function \ $F$ \ from \
$\kappa\times\kappa$ \ to \ $\kappa$ \ by stipulating that \ $F(\delta,\eta) =
(\chi(\eta))(\delta)$ \ if \ $\eta\in\kappa$ \ and \ $\delta\in 
\hbox{dom}(\chi(\eta)).$
\bigskip
Now fix \ $B\in I_{\kappa,\lambda}^+\cap P(A)$ \ and \ 
$\xi\in\kappa.$ \ Let \ $g\in
{}^\kappa\kappa$ \ be the increasing enumeration of the elements of the set \
$\{\cup(b\cap\kappa) : b\in B\}.$ \ Select \ $u\in 
P_\kappa(T)-\{\phi\}$ \ so that \
$g(\alpha)\leq\displaystyle\bigcup_{d\in u}(\psi_T(d))(\alpha)$ \ for 
all \ $\alpha<\kappa.$
\ Pick
\
$a\in B$ \ so that \ $\cup u\subseteq a$ \ and \ $\vmid{T\cap 
P(a)}>\xi.$ \ Then
$$g(\widetilde f_{T,a-\kappa}(\xi))<\cup(a\cap\kappa)$$
and
$$\xi = c_{f_{T,a-\kappa},\tau_{T,a-\kappa}}(g(\widetilde 
f_{T,a-\kappa}(\xi))) =
(\chi(\cup(a\cap\kappa))(g(\widetilde f_{T,a-\kappa}(\xi))) = F(g(\widetilde
f_{T,a-\kappa}(\xi)),\cup(a\cap\kappa)).$$
\hfill\carreblanc
\bigskip
{\bf Remark.} \ Theorems 14.2, 15.1 and 15.2 (as well as e.g. 
Theorems 9.2, 9.4, 12.1 and 12.2,
Propositions 9.3 and 11.3 and Corollary 13.2) are also true for \ 
$\kappa=\omega.$ \ This
gives (a) $\hbox{\nomgot d}\geq\hbox{\bf 
non}_\omega(J^+\hfltoutpetit{}{} [J^+]_\omega^2),$
\ and (b) if \ $\lambda\geq\hbox{\nomgot d},$ \ then  \
$I_{\omega,\lambda}^+\hfltoutpetit{\omega}{}\kern-3mm/\kern2mm
[I_{\omega,\lambda}^+]_\lambda^2$
\ and\ $I_{\omega,\lambda}^+ \hfltoutpetit{\omega}{\omega}\kern-3mm/\kern2mm
[I_{\omega,\lambda}^+]_\omega^2.$

\vskip 1,5cm

\vskip 1,5cm
\centerline{\Bbf{REFERENCES}}
\bigskip\bigskip
\newdimen\margeg \margeg=0pt
\def\bb#1&#2&#3&#4&#5&{\par{\parindent=0pt
     \advance\margeg by 1.1truecm\leftskip=\margeg
     {\everypar{\leftskip=\margeg}\smallbreak\noindent
     \hbox to 0pt{\hss\bf [#1]~~}{\bf #2 - }{\it#3~;} {#4.}\par\medskip
     #5 }
\medskip}}

\bb A&Y. ABE&Seminormal fine measures on $P_\kappa(\lambda)$& in : Proceedings
of the Sixth Asian Logic \break Conference, (C.T. Chong {\sl et al.,} 
eds.), World
Scientific, Singapore, 1998, pp. 1-12&&
\bb B&T. BARTOSZY\'NSKI&Combinatorial aspects of measure and
category&Fundamenta Mathema\-ticae 127 (1987), 225-239&&
\bb Bau1&J.E. BAUMGARTNER&Ineffability properties of cardinals I& in 
: Infinite and Finite
Sets (A. Hajnal, R. Rado and V.T. S\'os, eds.), Colloquia Mathematica 
Societatis J\'anos
Bolyai vol. 10, North-Holland, Amsterdam, 1975, pp. 109-130&&
\bb Bau2&J.E. BAUMGARTNER&Ineffability properties of cardinals II&in 
: Logic, Foundations
of\break Mathematics and Computability Theory (R.E. Butts and J. 
Hintikka, eds.), The
University of \break Western Ontario Series in Philosophy and Science 
vol. 9, Reidel,
Dordrecht (Holland), 1977, \break pp. 87-106&&
\bb BauT&J.E. BAUMGARTNER and A.D. TAYLOR&Partition theorems and 
ultrafilters&Transactions
of the American Mathematical Society 241 (1978), 283-309&&
\bb BauTW&J.E. BAUMGARTNER, A.D. TAYLOR and S. WAGON&Structural properties of
ideals&Dissertationes Mathematicae (Rozprawy Matematyczne) 197 (1982), 1-95&&
\bb Bl&A. BLASS&Ultrafilter mappings and their Dedekind 
cuts&Transactions of the American
Mathematical Society 188 (1974), 327-340&&
\bb CFMag&J. CUMMINGS, M. FOREMAN and M. MAGIDOR&Squares, scales and stationary
reflection&Journal of Mathematical Logic 1 (2001), 35-98&&
\bb CS&J. CUMMINGS and S. SHELAH&Cardinal invariants above the continuum&Annals
of Pure and Applied Logic 75 (1995), 251-268&&
\bb DMi&B. DUSHNIK and E.W. MILLER&Partially ordered sets&American Journal of
Mathematics 63 (1941), 600-610&&
\bb EHM\'aR&P. ERD\"OS, A. HAJNAL, A. M\'AT\'E and R. 
RADO&Combinatorial Set Theory :
Partition Relations for Cardinals&Studies in Logic and the Foundations of
Mathematics vol. 106, North-Holland, Amsterdam, 1984&&
\bb ER&P. ERD\"OS and R. RADO&A partition calculus in set 
theory&Bulletin of the
American Mathematical Society 62 (1956), 427-489&&
\bb Ka&A. KANAMORI&The Higher Infinite&Perspectives in Mathematical Logic,
Springer, Berlin, 1994&&
\bb K&K. KUNEN&Set Theory&North-Holland, Amsterdam, 1980&&
\bb Laf&C. LAFLAMME&Strong meager properties for filters&Fundamenta 
Mathematicae 146
(1995), 283-293&&
\bb L1&A. LANDVER&Singular Baire numbers and related topics&Ph. D. Thesis,
University of Wisconsin, Madison, Wisconsin, 1990&&
\bb L2&A. LANDVER&Baire numbers, uncountable Cohen sets and perfect-set
forcing&Journal of Symbolic Logic 57 (1992), 1086-1107&&
\bb M1&P. MATET&Combinatorics and forcing with distributive 
ideals&Annals of Pure and
Applied Logic 86 (1997), 137-201&&
\bb M2&P. MATET&The covering number for category and partition relations on
$P_\omega(\lambda)$&Fundamenta Mathematicae 171 (2002), 235-247&&
\bb M3&P. MATET&Partition relations for $\kappa$-normal ideals on
$P_\kappa(\lambda)$&Annals of Pure and Applied Logic 121 (2003), 89-111&&
\bb M4&P. MATET&A partition property of a mixed type for
$P_\kappa(\lambda)$&Mathematical Logic Quarterly 49 (2003), 615-628&&
\bb M5&P. MATET&Weak square bracket relations for 
$P_\kappa(\lambda)$&in preparation&&
\bb M6&P. MATET&Covering for category and combinatorics on 
$P_\kappa(\lambda)$&preprint&&
\bb MP1&P. MATET and J. PAWLIKOWSKI&Ideals over $\omega$ and cardinal 
invariants of
the continuum&Journal of Symbolic Logic 63 (1998), 1040-1054&&
\bb MP2&P. MATET and J. PAWLIKOWSKI&$Q$-pointness, $P$-pointness and 
feebleness of
ideals&Journal of Symbolic Logic 68 (2003), 235-261&&
\bb MP\'e&P. MATET and C. P\'EAN&Distributivity properties on
$P_\omega(\lambda)$&Discrete Mathematics, to appear&&
\bb MP\'eS1&P. MATET, C. P\'EAN and S. SHELAH&Cofinality of normal ideals on
$P_\kappa(\lambda)$ I&preprint&&
\bb MP\'eS2&P. MATET, C. P\'EAN and S. SHELAH&Cofinality of normal ideals on
$P_\kappa(\lambda)$ II&preprint&&
\bb MRoS&P. MATET, A. ROS\L ANOWSKI and S. SHELAH&Cofinality of the 
nonstationary
ideal&preprint&&
\bb Mil1&A.W. MILLER&A characterization of the least cardinal for 
which the Baire
category theorem fails&Proceedings of the American Mathematical 
Society 86 (1982),
498-502&&
\bb Mil2&A.W. MILLER&The Baire category theorem and cardinals of countable
cofinality&Journal of Symbolic Logic 47 (1982), 275-288&&
\bb S&S. SHELAH&On successors of singular cardinals&in : Logic 
Colloquium'78 (M. Boffa {\sl
et al.,} eds.), North-Holland, Amsterdam, 1979, pp. 357-380&&
\bb To1&S. TODORCEVIC&Partitioning pairs of countable  ordinals&Acta 
Mathematica 159
(1987), 261-294&&
\bb To2&S. TODORCEVIC&Coherent sequences&in : Handbook of Set Theory 
(M. Foreman, A.
Kanamori and M. Magidor, eds.), Kluwer, Dordrecht (Holland), to appear&&

\vskip 1,5cm
\settabs 3\columns
\+Universit\'e de Caen -CNRS&&The Hebrew University \cr
\+Math\'ematiques&&Institute of Mathematics \cr
\+BP 5186&&91904 Jerusalem \cr
\+14032 CAEN CEDEX&&Israel \cr
\+France&&\cr
\+e-mail : matet@math.unicaen.fr&&Rutgers University\cr
\+&&Departement of Mathematics\cr
\+&&New Brunswick, NJ 08854\cr
\+&&USA \cr
\+&&\cr
\+&&e-mail : shelah@math.huji.ac.il\cr

\end